\documentclass[10pt,journal,twocolumn]{IEEEtran}

\usepackage{amsmath,amssymb,epsf,epsfig,times}
\usepackage[all]{xy}
\usepackage{graphicx,color}
\usepackage{subfigure}
\usepackage{url}
\usepackage{cite}

\usepackage{keyval,times}
\usepackage[dvipdfm]{hyperref}
\hypersetup{%
   pdftitle={Cooperative Control},
   pdfauthor={Y. Cao and W. Ren},
   pdfkeywords={Cooperative Control, information consensus,
      multi-vehicle systems.},
   bookmarksnumbered,
   pdfstartview={FitH},
   colorlinks=true,
   breaklinks=true,
   citecolor=blue,
}%

\newtheorem{theorem}{Theorem}[section]
\newtheorem{lemma}{Lemma}[section]
\def\proof{\noindent{\it Proof: }}
\def\QED{\mbox{\rule[0pt]{1.5ex}{1.5ex}}}
\def\endproof{\hspace*{\fill}~\QED\par\endtrivlist\unskip}
\newcommand{\re}{\mathbb{R}}

\newcommand{\norm}[1]{\left\|#1\right\|}

\newcommand{\abs}[1]{\left|#1\right|}

\newcommand{\defeq}{\stackrel{\triangle}{=}}

\newtheorem{myremark}[theorem]{Remark}
\newtheorem{corollary}[theorem]{Corollary}

\newcommand{\Acal}{\mathcal{A}}
\newcommand{\Bcal}{\mathcal{B}}

\newcommand{\Lcal}{\mathcal{L}}

\newcommand{\Ncal}{\mathcal{N}}

\newcommand{\Gcal}{\mathcal{G}}

\newcommand{\onebf}{\mathbf{1}}
\newcommand{\zerobf}{\mathbf{0}}

\newcommand{\ith}{i{\text{th}}}

\newcommand{\OMIT}[1]{}

\begin{document}
\title{\LARGE \bf Finite-time Consensus for Multi-agent Networks with Unknown Inherent Nonlinear Dynamics}


\author{Yongcan Cao and Wei Ren 
\thanks{Y. Cao is with the Control Science Center of Excellence, Air Force Research Laboratory, Wright-Patterson AFB, OH 45433, USA. W. Ren is with the Department of Electrical Engineering,
        University of California, Riverside, CA 92521, USA. 
}
}


\maketitle

\maketitle 

\begin{abstract}
This paper focuses on analyzing the finite-time convergence of a nonlinear consensus algorithm for multi-agent networks with unknown inherent nonlinear dynamics. Due to the existence of the unknown inherent nonlinear dynamics, the stability analysis and the finite-time convergence analysis of the closed-loop system under the proposed consensus algorithm are more challenging than those under the well-studied consensus algorithms for known linear systems. For this purpose, we propose a novel stability tool based on a generalized comparison lemma. With the aid of the novel stability tool, it is shown that the proposed nonlinear consensus algorithm can guarantee finite-time convergence if the directed switching interaction graph has a directed spanning tree at each time interval. Specifically, the finite-time convergence is shown by comparing the closed-loop system under the proposed consensus algorithm with some well-designed closed-loop system whose stability properties are easier to obtain. Moreover, the stability and the finite-time convergence of the closed-loop system using the proposed consensus algorithm under a (general) directed switching interaction graph can even be guaranteed by the stability and the finite-time convergence of some special well-designed nonlinear closed-loop system under some special directed switching interaction graph, where each agent has at most one neighbor whose state is either the maximum of those states that are smaller than its own state or the minimum of those states that are larger than its own state. This provides a stimulating example for the potential applications of the proposed novel stability tool in the stability analysis of linear/nonlinear closed-loop systems by making use of known results in linear/nonlinear systems. For illustration of the theoretical result, we provide a simulation example.
\end{abstract}

\IEEEpeerreviewmaketitle

\section{Introduction} \label{sec:Introduction}

The past decade has witnessed an increasing research interest in the study of distributed cooperative control of multi-agent networks. The main objective is to design proper local controllers for a team of networked agents such that a desired group behavior can be accomplished. As one of the fundamental research topics in distributed cooperative control of multi-agent networks, consensus over multi-agent networks has been studied extensively. The main objective of consensus is to design distributed control algorithms such that a group of agents reach an agreement on some state of interest. Consensus has been investigated under various scenarios, including a deterministic interaction setting~\cite{SaberMurray04,RenBeard05_TAC,CaoMorseAnderson08,XiaoWang08}, a stochastic interaction setting~\cite{HatanoMesbahi05,Wu06,PorfiriStilwell07,SalehiJadbabaie08}, a sampled-data setting~\cite{XieLWJ09,XieLWJ09_2,GaoWang09,CaoRen10_IJRNC}, an asynchronous setting~\cite{CaoMorseAnderson08,FangAntsaklis08,XiaoWang08}, a quantization effect~\cite{KashyapBasarSrikant07,CarliBullo09}, and finite-time convergence~\cite{Cortes06,SundaramHadjicostis07,WangHong08_IFAC,HuiHaddadBhat08,XiaoWCG09,WangXiao10,WangHong10}.

Consensus with finite-time convergence, referred to as \emph{finite-time consensus}, means that consensus is achieved in finite time. In~\cite{Cortes06}, a nonsmooth consensus algorithm is proposed and the finite-time convergence of the closed-loop system is presented under an undirected fixed/switching interaction graph. In~\cite{HuiHaddadBhat08}, a continuous nonlinear consensus algorithm is proposed to guarantee the finite-time convergence under an undirected fixed interaction graph. In~\cite{WangXiao10}, the proposed algorithm in~\cite{HuiHaddadBhat08} is shown to guarantee finite-time consensus under an undirected switching interaction graph and a directed fixed interaction graph when each strongly connected component of the graph is detail-balanced. In~\cite{XiaoWCG09}, another continuous nonlinear consensus algorithm is proposed to guarantee the finite-time stability under a directed fixed interaction graph. In~\cite{WangHong10}, several nonlinear consensus algorithms are proposed to guarantee the finite-time $\chi$-consensus where the final equilibrium state can be predefined by designing the $\chi$ function. Note that in~\cite{Cortes06,HuiHaddadBhat08,XiaoWCG09,WangXiao10,WangHong10}, finite-time consensus is solved for single-integrator kinematics in a continuous-time setting. In~\cite{SundaramHadjicostis07}, finite-time consensus is studied for single-integrator kinematics in a discrete-time setting. To be specific, it is shown that the final equilibrium state can be computed after a finite number of time-steps. It is worth mentioning that finite-time consensus cannot be achieved when the continuous-time closed-loop systems under the finite-time consensus algorithms are discretized. In addition to the study of finite-time consensus for single-integrator kinematics as in~\cite{Cortes06,SundaramHadjicostis07,HuiHaddadBhat08,XiaoWCG09,WangXiao10,WangHong10}, a nonlinear distributed algorithm is proposed in~\cite{WangHong08_IFAC} to solve the finite-time consensus for double-integrator dynamics.

In the aforementioned papers, it is assumed that there is no inherent dynamics for the agents. However, inherent dynamics often exists for the agents in many practical systems. In the synchronization of complex dynamical networks~\cite{NishikawaMLH03,LuYuChen04,LuChen05,ZhouLuLu06}, the dynamics of each node in the complex networks is described by the sum of a continuously differentiable function describing the unknown inherent dynamics associated with the node and the coupling item identifying the connection between the node and other nodes. In~\cite{LuYuChen04}, chaos synchronization of general dynamical networks is studied under an undirected connected interaction graph, where the synchronization conditions are given in the form of matrix inequalities. In~\cite{LuChen05}, non-chaos synchronization of a time-varying complex dynamical network model is studied. To be specific, it is shown that synchronization is determined by the inner-coupling matrix as well as the eigenvalues and eigenvectors of the coupling configuration matrix characterizing the complex network. Similar to the results presented in~\cite{LuYuChen04}, the synchronization conditions in~\cite{LuChen05} are also given in the form of matrix inequalities. In~\cite{ZhouLuLu06}, adaptive algorithms are proposed to guarantee the local and global synchronization of an uncertain complex dynamic network. In~\cite{NishikawaMLH03}, the effect of heterogeneity in the synchronization of a complex network is investigated. It is shown that the ability of a scale-free network and a small-world network to synchronize decreases as the heterogeneity of the connectivity distribution increases. More details on the study of the synchronization of complex dynamical networks can be found in~\cite{Wu07,ChenWLL09}. Recently, the authors in~\cite{YuCCK09-} study second-order consensus of multi-agent systems with unknown Lipschitz nonlinear dynamics. Sufficient conditions are derived to guarantee second-order consensus under a directed fixed interaction graph. In~\cite{SuCWL10}, the authors propose a connectivity-preserving second-order consensus algorithm for multi-agent systems with unknown inherent nonlinear dynamics when there exists a virtual leader. In~\cite{YuChenCao11}, the authors study first-order consensus of multi-agent systems in the presence of unknown inherent nonlinear dynamics. Sufficient conditions are given to guarantee first-order consensus under a directed fixed interaction graph. Note that only asymptotical convergence is studied in the previous papers.

Motivated by the study of finite-time consensus for linear dynamics and asymptotical consensus with unknown inherent nonlinear dynamics, this paper studies finite-time consensus for multi-agent networks with unknown inherent nonlinear dynamics under a nonlinear consensus algorithm. Compared with the study of finite-time consensus for linear dynamics and asymptotic consensus with unknown inherent nonlinear dynamics, the study of finite-time consensus with unknown inherent nonlinear dynamics is more challenging because the two problems (\emph{i.e.}, finite-time consensus and asymptotical consensus with unknown inherent nonlinear dynamics) are considered simultaneously. In other words, the two problems are coupled as opposed to be separated. It is worth noting that: (1) finite-time consensus with unknown inherent nonlinear dynamics has not been investigated in the existing literature; and (2) the approaches used in the study of finite-time consensus for linear dynamics and asymptotic consensus with unknown inherent nonlinear dynamics are not applicable to the study of finite-time consensus for multi-agent networks with unknown inherent nonlinear dynamics. Therefore, it is potentially challenging to conduct the stability analysis and the finite-time convergence analysis of the closed-loop systems under the proposed nonlinear consensus algorithm. To facilitate the stability analysis and the finite-time convergence analysis, we propose a novel stability tool based on a generalized comparison lemma. By using the novel stability tool, we show that the proposed nonlinear consensus algorithm can guarantee finite-time convergence for multi-agent networks with unknown inherent nonlinear dynamics if the directed switching interaction graph has a directed spanning tree at each time interval. Specifically, the finite-time convergence is shown by comparing the closed-loop system under the proposed consensus algorithm with some well-designed closed-loop system whose stability properties are easier to obtain. Moreover, by using the novel stability tool, the stability and the finite-time convergence of the closed-loop system using the proposed consensus algorithm under a (general) directed switching interaction graph can even be guaranteed by the stability and the finite-time convergence of some well-designed nonlinear closed-loop system under some special directed switching interaction graph, where each agent has at most one neighbor whose state is either the maximum of those states that are smaller than its own state or the minimum of those states that are larger than its own state. This provides a stimulating example for the potential applications of the proposed novel stability tool. As a byproduct, in the absence of the unknown inherent nonlinear dynamics, the proposed nonlinear consensus algorithm can still guarantee finite-time convergence if the directed switching interaction graph has a directed spanning tree at each time interval. This extends the existing research on the study of finite-time consensus for single-integrator kinematics to a more general case where a milder condition on the interaction graph is required.


The remainder of this paper is organized as follows. In Section~\ref{sec: pre}, we briefly review notations used in this paper, the graph theory notions, and the problem to be solved. In Section~\ref{sec:General_comparison_lemma}, we propose a novel stability tool based on a generalized
comparison lemma and present the basic steps to use the novel stability tool. In Section~\ref{sec:directed}, we analyze the stability of the proposed consensus algorithm by using the novel stability tool under a directed switching interaction graph. Then a simulation example is given in Section~\ref{sec:sim} to further validate the result obtained in Section~\ref{sec:directed}. Finally, Section~\ref{sec:conclusion} is given to summarize the contribution of the paper.

\section{Preliminaries and Problem Statement}\label{sec: pre}
\subsection{Notations}
We use $\re$ to denote the set of real numbers. $\zerobf_n\in\re^n$ is used to denote the $n\times 1$ all-zero column vector. $\onebf_n\in\re^n$ is the $n\times1$ all-one column vector. $I_n\in\re^{n\times n}$ is used to denote the identity matrix. $\norm{\cdot}$ is used to denote the $2$-norm of a vector. Define $\text{sig}(x)^\alpha\defeq\text{sgn}(x)\abs{x}^\alpha$, where $\text{sgn}(\cdot)$ is used to denote the signum function. Note that $\text{sig}(x)^\alpha$ is continuous with respect to $x$ when $\alpha>0$. Let $f:[0,\infty)\mapsto J\subseteq\re^n$ be a continuous function. The upper Dini derivative of $f(t)$ is given by $D^+f(t)=\limsup_{h\to 0+}\frac{1}{h}[f(t+h)-f(t)]$. Given a function $f(x(t))$, the upper Dini derivative of $f(x(t))$ is defined as $D^+f(x(t))=\limsup_{h\to 0+}\frac{1}{h}[f(x(t+h))-f(x(t))]$. A function $f(t,x)$ is \emph{locally $\chi$-Lipschitz} (also known as locally H\"{o}lder continuous) in $x$ if there exist real positive constants $C$, $\chi$, and $\epsilon$ such that $\norm{f(t,x)-f(t,y)}\leq C\norm{x-y}^\chi$ for all $y\in \Bcal(x,\epsilon)$, where $\Bcal(x,\epsilon)$ denotes the ball centered at $x$ with a radius $\epsilon$.


\subsection{Graph Theory Notions}
For a team of $n$ agents, the interaction among them can be modeled by a directed graph $\mathcal {G}=(\mathcal {V,W})$, where $\mathcal{V}=\{1,2,\ldots,n\}$ and $\mathcal{W}\subseteq \mathcal{V}^{2}$ represent, respectively, the agent set and the edge set. An edge denoted as $(i, j)$ means that agent $j$ can obtain information from agent $i$. That is, agent $i$ is a neighbor of agent $j$. We use $\Ncal_j$ to denote the neighbor set of agent $j$. A directed path is a sequence of edges of the form $(i_1,i_2),(i_2,i_3),\ldots,$ where $i_k\in\mathcal{V},~k=1,2,\cdots$. A directed graph has a directed spanning tree if there exists at least one agent that has directed paths to all other agents.


Two matrices are frequently used to represent the interaction graph: the adjacency matrix $\Acal=[a_{ij}]\in \re^{n\times n}$ with
$a_{ij}>0$ if $(j,i)\in \mathcal{W}$ and $a_{ij}=0$ otherwise, and the Laplacian matrix $\Lcal=[\ell_{ij}]\in
\re^{n\times n}$ with $\ell_{ii}=\sum_{j=1,j\neq i}^na_{ij}$ and $\ell_{ij}=-a_{ij}$, $i\neq j$. In particular, we let $a_{ii}=0,~i=1,\ldots,n,$ (\textit{i.e.}, agent $i$ is not a neighbor of itself). It is straightforward to verify that $\Lcal$ has at least one eigenvalue equal to $0$ with a corresponding right eigenvector $\onebf_n$.


\subsection{Problem Statement}
Consider a group of $n$ agents with dynamics given by
\begin{equation}\label{eq:kinmatics}
  \dot{r}_i=\phi(t,r_i)+u_i,\quad i=1,\ldots,n,
\end{equation}
where $r_i\in\re^m$ is the state of the $\ith$ agent, $\phi(t,r_i)\in\re^m$ is the unknown inherent Lipschitz nonlinear dynamics for the $\ith$ agent, and
$u_i\in\re^m$ is the control input for the $\ith$ agent. Note that the nonlinear term
$\phi(t,r_i)\defeq [\phi(t,r_i^{(1)}),\ldots,\phi(t,r_i^{(m)})]$ is unknown, but is assumed to satisfy
\begin{equation}\label{eq:lipshitz}
\abs{\phi(t,r_i^{(k)})-\phi(t,r_j^{(k)})}\leq \gamma\abs{r_i^{(k)}-r_j^{(k)}},\quad k=1,\ldots,m,
\end{equation}
where $r_i^{(k)}$ is the $k$th component of $r_i$ and $\gamma$ is a known positive constant.
The objective is to design $u_i$ such that $\norm{r_i(t)-r_j(t)}\to 0$ in finite time for all $i,j=1,\ldots,n$. That is, all agents' states reach consensus in finite time. Due to the existence of the unknown nonlinear term $\phi(t,r_i)$ in~\eqref{eq:kinmatics}, the consensus value, in general, is not constant, which is different from the case when the unknown inherent nonlinear dynamics do not exist. Moreover, since $\phi(t,r_i)$ is unknown, it is impossible to introduce the term $-\phi(t,r_i)$ in the control algorithm to cancel the unknown nonlinear term $\phi(t,r_i)$.

It is well known that linear algorithms normally cannot guarantee finite-time convergence. We propose the following \emph{nonlinear} finite-time consensus algorithm for~\eqref{eq:kinmatics} as
\begin{equation}
  u_i = -\beta\sum_{j=1}^n a_{ij}(t) \text{sig}(r_i-r_j)^{\alpha(\norm{r_i-r_j})},\label{eq:control-new}
\end{equation}
where $\beta$ is a positive constant, $a_{ij}(t)$ is the $(i,j)$th entry of the adjacency matrix $\Acal(t)$ associated with the graph $\Gcal(t)$ characterizing the interaction among the $n$ agents at time $t$, $\text{sig}(\cdot)$ is defined componentwise, and $\alpha(\norm{r_i-r_j})$ is defined such that
\begin{align}\label{eq:alpha-func}
\alpha(\norm{r_i-r_j})=\left\{
\begin{array} {ll}
\alpha^\star,&0\leq\norm{r_i-r_j}<1,\\
1,&\norm{r_i-r_j}\geq 1,
\end{array}\right.
\end{align}
where $\alpha^\star\in(0,1)$ is a positive constant. The main idea behind~\eqref{eq:control-new} is that an agent uses the information $r_i-r_j$ for $j\in\Ncal_i(t)$ when $\norm{r_i-r_j}\geq 1$ while uses $\text{sig}(r_i-r_j)^{\alpha^\star}$ for $j\in\Ncal_i(t)$ when $\norm{r_i-r_j}< 1$. The objective of using $r_i-r_j$ is for the stabilization and the objective of using $\text{sig}(r_i-r_j)^{\alpha^\star}$ is to further guarantee finite-time convergence.

In this paper, we assume that the adjacency matrix $\Acal(t)$ is constant for $t\in[t_i,t_{i+1})$ and switches at time
$t_{i+1}$, $i=0,1,\ldots$, where $t_0=0$. Let $\Gcal_i$, $\Acal_i$, and $\Lcal_i$ denote, respectively, the
directed graph, the adjacency matrix, and the Laplacian matrix associated with the $n$ agents for
$t\in[t_i,t_{i+1})$. We assume that $t_{i+1}-t_i\geq t_L$, where $t_L$ is a positive constant. We also assume that each nonzero and hence positive entry of $\Acal_i$ has a lower bound $\underline{a}$ and an upper bound $\overline{a}$, where $\underline{a}$ and $\overline{a}$ are positive constants.

Note that $\alpha(\norm{r_i-r_j})$ is discontinuous with respect to $\norm{r_i-r_j}$ when $\norm{r_i-r_j}=1$. However, $\text{sig}(r_i-r_j)^{\alpha(\norm{r_i-r_j})}$ is continuous with respect to $\norm{r_i-r_j}$ when $\norm{r_i-r_j}=1$ by recalling the definition of $\text{sig}(\cdot)^\alpha$. Accordingly, the right-hand side of~\eqref{eq:kinmatics} using~\eqref{eq:control-new} is discontinuous when the interaction graph switches. Because the interaction graph switches only at some distinct time instants, the set of the discontinuity points for the right-hand side of~\eqref{eq:kinmatics} using~\eqref{eq:control-new} has measure zero. It follows from~\cite{Filippov88} that the Caratheodory solution to the closed-loop system of~\eqref{eq:kinmatics} using~\eqref{eq:control-new} exists for an arbitrary initial condition. In addition, the solution is an absolutely continuous function that satisfy~\eqref{eq:kinmatics} using~\eqref{eq:control-new} for almost all $t$. In the following of this paper, we consider the solution of~\eqref{eq:kinmatics} using~\eqref{eq:control-new} in the Caratheodory sense.


\begin{myremark}
A special case of the nonlinear algorithm~\eqref{eq:control-new} with $\alpha(\norm{r_i-r_j})\equiv \alpha^\star$, where $\alpha^\star\in(0,1)$ is a constant scalar, is used in~\cite{HuiHaddadBhat08} to solve finite-time consensus for single-integrator agents without the unknown inherent nonlinear dynamics under an undirected fixed interaction graph. A variant of the nonlinear algorithm~\eqref{eq:control-new} is used in~\cite{WangXiao10} to solve finite-time consensus for single-integrator agents without the unknown inherent nonlinear dynamics under an undirected switching interaction graph and a directed fixed interaction graph when each strongly connected component of the graph is detailed-balanced. Different from~\cite{HuiHaddadBhat08,WangXiao10}, we investigate finite-time consensus for agents with the unknown inherent nonlinear dynamics \emph{under a general directed switching interaction graph}. The approaches used in~\cite{HuiHaddadBhat08,WangXiao10} to conduct the stability analysis are not applicable to the problem studied in this paper.
\end{myremark}

\section{A Novel Stability Tool By Comparison}\label{sec:General_comparison_lemma}
In this section, we propose a novel stability tool based on a generalized comparison lemma, which is used in the following of the paper for the stability analysis of the proposed algorithm~\eqref{eq:control-new}.

\begin{theorem} \label{th:comparison-vector-general}
Consider a nonlinear system given by
\begin{equation}\label{eq:sys2}
\dot{x}=g(t,x),\quad x\in\re^m,
\end{equation}
where $g(t,x)$ is locally $\chi$-Lipschitz in $x$ and is continuous in $t$. Let $G(x):\re^m\mapsto \re$ be a nonnegative function satisfying that $G(x)=0$ for any $t\geq \underline{t}$ if $G(x)=0$ at $t=\underline{t}$.\footnote{$G(x)$ is a function of $t$ since $x$ is a function of $t$.} Then $G(x)\to 0$ in finite time (respectively, as $t\to\infty$) for any initial state if there exist a nonnegative function $F(z):\re^m\mapsto \re$ and a function $f(t,z)$ that is locally $\chi$-Lipschitz in $z$ and is continuous in $t$ such that the following conditions are satisfied:
\begin{itemize}
  \item[1.] $G(\xi)\leq F(\xi)$ for any $\xi\in\re^m$;
  \item[2.] If $F(\xi)\neq 0$ (and hence $>0$), $\frac{\text{d}G(\xi)}{\text{d}\xi}g(t,\xi)<\frac{\text{d}F(\xi)}{\text{d}\xi}f(t,\xi)$ holds;
  \item[3.] $\frac{\text{d}F(z)}{\text{d}z}f(t,z)$ is continuous with respect to $z$;
  \item[4.] At least one solution exists for
  \begin{align}\label{eq:equality-derivative}
      \frac{\text{d}F(z)}{\text{d}z}\dot{z}=\frac{\text{d}F(z)}{\text{d}z}f(t,z).
  \end{align}
  For any $z$ satisfying~\eqref{eq:equality-derivative}, $F(z)\to 0$ in finite time (respectively, as $t\to\infty$) for an arbitrary initial state.
  \item[5.] When $F(z(0))\leq F_0$, either of the following two cases holds:
   \begin{itemize}
   \item[(1).] $z(0)$ is bounded for any positive constant $F_0$;
   \item[(2).] $\max_i z_i(0)-\min_i z_i(0)$ is bounded for any positive constant $F_0$ and $F(z)$ satisfying~\eqref{eq:equality-derivative} is invariant when $z(0)=z(0)+\eta\textbf{1}_m$ for any positive constant $\eta$, where $z_i$ is the $i$th entry of $z$;
   \end{itemize}
\end{itemize}
\end{theorem}
\proof 
When $g(t,x)$ is locally $\chi$-Lipschitz in $x$ and is continuous in $t$, we have that~\eqref{eq:sys2} has a unique solution. The main idea of the proof is to show that there exist time intervals $[0,t_1),[t_1,t_2),\ldots$ such that during any time interval $[t_i,t_{i+1}),~i=0,1,\ldots,$ where $t_0=0$, there exists a proper initial state $\hat{z}^i(0)$ satisfying $F(\hat{z}^i(0))\leq F(x(0))$ such that
$G(x)\leq \max_{z\in \Pi} F(z)$
for all $t\in[t_i,t_{i+1})$, where $\Pi\defeq \{\pi|\frac{\text{d}F(\pi)}{\text{d}\pi}\dot{\pi}=\frac{\text{d}F(\pi)}{\text{d}\pi}f(t,\pi)\}$.\footnote{Note that both $G(x)$ and $F(z)$ are functions of $t$ because $x$ and $z$ are functions of $t$.}

First, we show that there exists $t_1$ such that there exists a proper initial state $\hat{z}^0(0)$ satisfying $F(\hat{z}^0(0))\leq F(x(0))$ to guarantee that $G(x)\leq \max_{z\in \Pi} F(z)$ for all $t\in[0,t_1)$. Because $F(z)$ is differentiable with respect to $z$ and $z$ is differentiable with respect to $t$, $F(z)$ is continuous with respect to $z$ and $z$ is continuous with respect to $t$. Therefore, $F(z)$ is continuous with respect to $t$. Then there exists $\tau_1$ such that $\max_{z\in\Pi} F(z)$ is well defined for $t\in[0,\tau_1)$. In other words, there exists a solution $z^\star$ such that $F(z^\star)\geq F(z)$ for any $z\in\Pi$ and any $t\in[0,\tau_1)$. Let the initial state be chosen as $\hat{z}^0(0)=x(0)$. From Condition 1, it is apparent that $G(x(0))\leq F(x(0))=F(\hat{z}^0(0))$. When $F(\hat{z}^0(0))=0$, it follows from Condition 1 that $G(x(0))=0$. Because $G(x)=0$ for any $t\geq \underline{t}$ if $G(x)=0$ at $t=\underline{t}$, it follows that $G(x)=0$ for $t\geq 0$. Then the theorem holds trivially. We next consider the case when $F(\hat{z}^0(0))\neq0$. From Condition 2, it follows that
\begin{align*}
&\frac{\text{d}G(x)}{\text{d}t}=\frac{\text{d}G(x)}{\text{d}x}g(t,x)\\
<&\frac{\text{d}[\max_{x\in \Pi}F(x)]}{\text{d}x}f(t,x)=\frac{\text{d}[\max_{x\in \Pi}F(x)]}{\text{d}t}
\end{align*}
at $t=0$ when $F(x(0))\neq 0$. Then there exists $T_0>0$ such that $G(x)\leq \max_{z\in \Pi}F(z)$ for $t\in [0,T_0)\subseteq[0,\tau_1)$. Noting that the relationship between $x(T_0)$ and $z(T_0)$ is unknown, we might not choose $T_0$ as $t_1$ because the relationship between $\frac{\text{d}G(x)}{\text{d}t}$ and $\frac{\text{d}[\max_{x\in \Pi}F(x)]}{\text{d}t}$ might be unknown under the condition of the theorem. Accordingly, the existence of $t_2$ is not guaranteed. Here, the objective is to choose a proper initial state $\hat{z}^1(0)$ and $t_1$ such that $F(\hat{z}^1(0))\leq F(x(0))$ and $x(t_1)=z(t_1)$, where $t_1\in(0,T_0)$. We next show that there exist $t_1\in [0,T_0)$ and $\hat{z}^1(0)$ such that $F(\hat{z}^1(0))\leq F(x(0))$ and $z(t_1)=x(t_1)$. Denote
$
\eta^0(t)\defeq \frac{\text{d}[\max_{z\in \Pi}F(z)]}{\text{d}z}f(t,z)|_{z(0)=\hat{z}^0(0)}-\frac{\text{d}G(x)}{\text{d}x}g(t,x)|_{x(0)=\hat{z}^0(0)}
$
and
$
\Delta^0_F(t)\defeq \frac{\text{d}[\max_{z\in \Pi}F(z)]}{\text{d}z}f(t,z)|_{z(0)=x(0)}-\frac{\text{d}[\max_{z\in \Pi}F(z)]}{\text{d}z}f(t,z)|_{z(t_1)=x(t_1)}.
$
Because $\eta^0(t)>0$ when $t=0$, $\Delta^0_F(t)=0$ when $t_1=0$, $\eta^0(t)$ is continuous with respect to $t$, and $\Delta^0_F(t)$ is continuous with respect to $z$ (Condition 3), there exists $t_1$ such that $\Delta^0_F(t)\leq \eta^0(t)$ when $t\in [0,t_1)\subseteq[0,T_0)$. Therefore, $F(z(0))|_{z(t_1)=x(t_1)}\leq F(z(0))|_{z(0)=x(0)}$. The initial state that guarantees $z(t_1)=x(t_1)$ is designated as $\hat{z}^1(0)$.

Next, we show that there exists $t_2$ such that $G(x)\leq \max_{z\in \Pi}F(z)$ for all $t\in[t_1,t_2)$ under the initial state $\hat{z}^1(0)$. We only consider the case when $F(z(t_1))\neq0$ since the theorem holds trivially if $F(z(t_1))=0$ by noting that $0\leq G(x(t_1))\leq F(z(t_1))=0$ and $G(x)=0$ for any $t\geq \underline{t}$ if $G(x)=0$ at $t=\underline{t}$. By following the previous analysis, $\frac{\text{d}G(x)}{\text{d}t}<\frac{\text{d}[\max_{z\in \Pi}F(z)]}{\text{d}t}$ at $t=t_1$. Then there exists $T_1>0$ such that $G(x)\leq \max_{z\in \Pi}F(z)$ for $t\in [t_1,t_1+T_1)$. Again, we next show that there exists a proper initial state $\hat{z}^2(0)$ and $t_2$ such that $F(\hat{z}^2(0))\leq F(x(0))$ and $x(t_2)=z(t_2)$, where $t_2\in(t_1,t_1+T_1)$. Denote
$
\eta^1(t)\defeq \frac{\text{d}[\max_{z\in \Pi}F(z)]}{\text{d}z}f(t,z)|_{z(0)=\hat{z}^1(0)}-\frac{\text{d}G(x)}{\text{d}x}g(t,x)|_{x(0)=\hat{z}^0(0)}
$
and
$
\Delta^1_F(t)\defeq \frac{\text{d}[\max_{z\in \Pi}F(z)]}{\text{d}z}f(t,z)|_{z(t_1)=x(t_1)}-\frac{\text{d}[\max_{z\in \Pi}F(z)]}{\text{d}z}f(t,z)|_{z(t_2)=x(t_2)}.
$
Because $\eta^1(t)>0$ when $t=t_1$, $\Delta^1_F(t)=0$ when $t_2=t_1$, $\eta^1(t)$ is continuous with respect to $t$, and $\Delta^1_F(t)$ is continuous with respect to $z$, there exists $t_2$ such that $\Delta^1_F(t)\leq \eta^1(t)$ when $t\in [0,t_1+t_2)\subseteq[0,T_2)$. Therefore, $F(z(0))|_{z(t_2)=x(t_2)}\leq F(z(0))|_{z(0)=\hat{z}^1(0)}\leq F(0,z(0))|_{z(0)=x(0)}$. The initial state that guarantees $z(t_2)=x(t_2)$ is designated as $\hat{z}^2(0)$.

By following a similar analysis, it can be shown that whenever $G(x)\neq 0$ at some $t_i$, there always exists $t_{i+1}>t_i$ such that there exists a proper initial state $\hat{z}^i(0)$ satisfying $F(\hat{z}^i(0))\leq F(x(0))$ that can guarantee $G(x)\leq \max_{z\in \Pi}F(z)$ for $t\in[t_i,t_{i+1})$. When the Subcondition (1) under Condition 5 is satisfied, $F(z(0))\leq F(x(0))$ implies that $z(0)$ is always from a bounded set. Combining with Condition 4 shows that $F(z)\to 0$ in finite time (respectively, as $t\to\infty$) for any $z\in\Pi$, which implies that $G(x)\to 0$ in finite time (respectively, as $t\to\infty$). When the Subcondition (2) under Condition 5 is satisfied, if there exists an initial condition $z(0)$ under which a desired trajectory $F(z)$ can be achieved, the same desired trajectory can be achieved if adding each $z_i(0)$ by a common constant. Therefore, there always exists a bounded set of $z(0)$ such that any trajectory $F(z)$ with the initial state satisfying $F(z(0))\leq F(x(0))$ can be achieved when the corresponding initial condition is chosen properly from the set. Again, $G(x)\to 0$ in finite time (respectively, as $t\to\infty$) by considering Condition 4.

This completes the proof of the theorem.
\endproof

\begin{myremark}
There are five conditions in Theorem~\ref{th:comparison-vector-general}. Each condition has its own unique purpose. Condition 1 is to guarantee that $G(x)-F(z)$ are nonpositive once $x$ and $z$ are equivalent. Condition 2 is to guarantee that the derivative of $G(x)$ along the solution of~\eqref{eq:sys2} is smaller than the derivative of $F(z)$ along the solution of~\eqref{eq:equality-derivative} once $x=z$ and $F(z)$ is positive. Conditions 1 and 2 guarantee that once $x=z$ and $F(z)$ is positive, there exists an time interval during which $G(x)-F(z)\leq 0$. However, Conditions 1 and 2 might not guarantee the existence of a proper initial state $z(0)$ such that $F(z(0))|_{z(t_i)=x(t_i)}\leq F(z(0))|_{z(0)=x(0)}$ for some future time instant $t_i$ unless Condition 3 is imposed. Condition 4 guarantees that $F(z)\to 0$ in finite time (respectively, as $t\to\infty$) for any given initial state. Condition 5 guarantees that for any positive $F_0$, the solutions to~\eqref{eq:equality-derivative} with any initial state $z(0)$ satisfying $F(z(0))\leq F_0$ are always within the set of solutions to~\eqref{eq:equality-derivative} with the initial state $z(0)$ being chosen from some fixed bounded set. Without Condition 5, it is possible that the fixed bounded set becomes a dynamic unbounded set. Consequently, $F(z(0))$ might become unbounded. Then $F(z(t))$ might become unbounded as well. As such, Conditions 4 and 5 guarantee that $F(z)$ are bounded and $F(z)\to 0$ in finite time (respectively, as $t\to\infty$) for every $z(0)$ satisfying $F(z(0))\leq F_0$. Therefore, the five conditions are required in Theorem~\ref{th:comparison-vector-general}. Similarly, the five conditions are required in the following Corollary~\ref{co:comparison-vector-general} as well.\endproof
\end{myremark}

Although $F(z)$ and $G(x)$ are essentially functions of $t$ by noting that both $z$ and $x$ are functions of $t$, we choose to use $F(z)$ and $G(x)$ rather than $F(t)$ and $G(t)$ because the omission of $z$ and $x$ could substantially affect the demonstration of the close connection between $\frac{\text{d}F(z)}{\text{d}t}$ and $\frac{\text{d}G(x)}{\text{d}t}$. The close connection is primarily given in Conditions 1 and 2 in Theorem~\ref{th:comparison-vector-general}. Condition 1 means that $F(\cdot)-G(\cdot)$ is nonnegative when a common state is applied to both $F$ and $G$. Condition 2 means that $\frac{\text{d}}{\text{d}t}[F(\cdot)-G(\cdot)]$ is positive when a common state is applied to both $F$ and $G$ if $F$ is positive. Nevertheless, the existence of common states for $F$ and $G$ is neither necessary nor feasible for two physical systems satisfying the conditions in Theorem~\ref{th:comparison-vector-general}. In other words, Conditions 1 and 2 reflect only the \emph{structural} requirement (what structural properties the system have when the states are the same). In order to better demonstrate the structural relationship between $\frac{\text{d}F(z)}{\text{d}t}$ and $\frac{\text{d}G(x)}{\text{d}t}$, the notation $F(z)$ and $G(x)$ (or the alike) will be employed in the following part of the paper.

Compared with the comparison lemmas for vector differential equations in~\cite{MichelMiller77,Siljak78} (see Theorem 2.6.11 in~\cite{MichelMiller77} and Theorem 2.1 in~\cite{Siljak78}), the generalized comparison lemma in Theorem~\ref{th:comparison-vector-general} does not require the \textit{quasimonotone} condition (see Definition 2.6.9 in~\cite{MichelMiller77} and Definition 2.9 in~\cite{Siljak78}). In addition, we use scalar functions in Theorem~\ref{th:comparison-vector-general} but vector functions are used in~\cite{MichelMiller77,Siljak78}.

Note that Theorem~\ref{th:comparison-vector-general} considers the case when the right-hand side of~\eqref{eq:sys2} is smooth and $F(z)$ and $G(x)$ are differentiable with respect to $t$. We next present an extension of Theorem~\ref{th:comparison-vector-general} by considering the case when the right-hand side of~\eqref{eq:sys2} are smooth almost everywhere and $F(z)$ and $G(x)$ are differentiable with respect to $t$ almost everywhere. In this case, the solution to~\eqref{eq:sys2} is considered in the Caratheodory sense. The Caratheodory solution of~\eqref{eq:sys2} exists because we will consider the case where the set of the discontinuity points for the right-hand side of~\eqref{eq:sys2} has measure zero. This relaxation has a similar counterpart in~\cite{MichelMiller77} (c.f. Theorem 2.6.11).

\begin{corollary} \label{co:comparison-vector-general}
Consider a nonlinear system given by~\eqref{eq:sys2}
where $g(t,x)$ is locally $\chi$-Lipschitz in $x$ almost everywhere and is continuous in $t$. Let $G(x):\re^m\mapsto \re$ be a nonnegative function satisfying that $G(x)=0$ for any $t\geq \underline{t}$ if $G(x)=0$ at $t=\underline{t}$. Then $G(x)\to 0$ in finite time (respectively, as $t\to\infty$) for any initial state if there exist a nonnegative function $F(z):\re^m\mapsto \re$ and a function $f(t,z)$ that is locally $\chi$-Lipschitz in $z$ almost everywhere and is continuous in $t$ such that the following conditions are satisfied:
\newline
1. $G(\xi)\leq F(\xi)$ for any $\xi\in\re^m$;\newline
2. If $F(\xi(t))\neq 0$ (and hence $>0$) and $z(t)=x(t)=\xi(t)$,
$\limsup_{h\to 0+}\frac{G(x(t+h))-G(x(t))}{h}|_{\dot{x}(t)=g(t,x(t))}<\limsup_{h\to 0+}\frac{F(z(t+h))-F(z(t))}{h}|_{\dot{z}(t)=f(t,z(t))}$
holds almost everywhere;\newline
3. $\limsup_{h\to 0+}\frac{F(z(t+h))-F(z(t))}{h}|_{\dot{z}(t)=f(t,z(t))}$ is continuous with respect to $z$ almost everywhere;\newline
4. At least one (Caratheodory) solution $z$ exists for
\begin{align}\label{eq:equality-Dini}
&\limsup_{h\to 0+}\frac{F(z(t+h))-F(z(t))}{h}\notag\\
=&\limsup_{h\to 0+}\frac{F(\xi(t+h))-F(\xi(t))}{h}|_{\dot{\xi}(t)=f(t,\xi(t))}.
\end{align}
Furthermore, for any (Caratheodory) solution $z$ satisfying~\eqref{eq:equality-Dini}, $F(z)\to 0$ in finite time (respectively, as $t\to\infty$) for an arbitrary initial state.\newline
5. When $F(z(0))\leq F_0$, either of the following two cases holds:
   \begin{itemize}
   \item[(1).] $z(0)$ is bounded for any positive constant $F_0$;
   \item[(2).] $\max_i z_i(0)-\min_i z_i(0)$ is bounded for any positive constant $F_0$ and $F(z)$ satisfying~\eqref{eq:equality-Dini} is invariant when $z(0)=z(0)+\eta\textbf{1}_m$ for any positive constant $\eta$, where $z_i$ is the $i$th entry of $z$;
   \end{itemize}
\end{corollary}
\proof It can be observed that Conditions 1 and 5 are exactly the same as those in Theorem~\ref{th:comparison-vector-general}, while Conditions
2, 3, and 4 are the same as those in Theorem~\ref{th:comparison-vector-general} except that the use of derivative in Conditions
2, 3, and 4 is replaced by the
use of upper Dini derivative and the statements in Conditions 2 and 3 are assumed to be satisfied almost everywhere instead of everywhere.
Since the exclusion of the discontinuity points with measure zero does not change the solutions of~\eqref{eq:equality-Dini} in the Caratheodory sense, the proof of the corollary is the same as that of Theorem~\ref{th:comparison-vector-general} when the discontinuity points with measure zero are excluded.
\endproof

Theorem~\ref{th:comparison-vector-general} and Corollary~\ref{co:comparison-vector-general} reflect the connection between the stability of one system and the stability of other system(s) when certain conditions, as mentioned in Theorem~\ref{th:comparison-vector-general} and Corollary~\ref{co:comparison-vector-general}, are satisfied. The basic idea is to show that some nonnegative function will approach zero in finite time (respectively, as $t\to\infty$), which has a similar motivation to that of the Lyapunov-based stability analysis. However, substantial differences can be found between them. Recall that the asymptotic/finite-time stability of a closed-loop system can be guaranteed via the Lyapunov-based stability analysis when there exists a proper Lyapunov function with a negative-semidefinite (and the alike) derivative. In contrast, it is \emph{not} required in Theorem~\ref{th:comparison-vector-general} and Corollary~\ref{co:comparison-vector-general} that some nonnegative function (\textit{i.e.,} $G(\omega)$ in Theorem~\ref{th:comparison-vector-general} and Corollary~\ref{co:comparison-vector-general}) has a negative-semidefinite (and the alike) derivative. Instead, the derivative of the nonnegative function is studied through a proper comparison with that of some nonnegative function under some properly designed systems, whose stability can be obtained or is already known. More precisely and concisely, Theorem~\ref{th:comparison-vector-general} and Corollary~\ref{co:comparison-vector-general} can be used to analyze the stability of some closed-loop system given that the conditions in Theorem~\ref{th:comparison-vector-general} and Corollary~\ref{co:comparison-vector-general} are satisfied even if it is difficult to find a Lyapunov function for the system or a Lyapunov function for the system does not exist (since the existence of a Lyapunov function for a stable system is not always guaranteed). This is the main contribution of the proposed novel stability tool based on the generalized comparison lemma. Before applying Theorem~\ref{th:comparison-vector-general} and Corollary~\ref{co:comparison-vector-general} in the stability analysis of~\eqref{eq:kinmatics} using the proposed consensus algorithm~\eqref{eq:control-new}, we next present the major steps for the use of Theorem~\ref{th:comparison-vector-general} and Corollary~\ref{co:comparison-vector-general}.

\textbf{Major steps to use Theorem~\ref{th:comparison-vector-general} and Corollary~\ref{co:comparison-vector-general}:}
\begin{itemize}
\item[1.] Choose a proper nonnegative function $G(x)$ based on the closed-loop system dynamics and compute $\limsup_{h\to 0+}\frac{G(x(t+h))-G(x(t))}{h}|_{\dot{x}(t)=g(t,x(t))}$.\footnote{Note that $\limsup_{h\to 0+}\frac{G(x(t+h))-G(x(t))}{h}|_{\dot{x}(t)=g(t,x(t))}=\frac{\text{d}G(\xi)}{\text{d}\xi}g(t,\xi)$ when $G(x)$ is differentiable.} Similar to the idea behind the Lyapunov functions, it is normally a good idea to choose a quadratic function $x^TQx$ with $Q$ being symmetric positive semi-definite, $\max_ix_i-\min_i x_i$ with $x_i$ being the $\ith$ entry of $x$, $\max_i\abs{x_i}$, and etc.
\item[2.] Design a proper function $F(z)$ satisfying Condition 1 and the corresponding closed-loop system(s) satisfying~\eqref{eq:equality-derivative} or~\eqref{eq:equality-Dini} such that Conditions 2-5 are satisfied.
\end{itemize}
It seems that the two steps are quite independent from the previous descriptions. However, they are closely related because the choice of $G(x)$ in Step 1 has a direct impact on whether the design of a proper function and the corresponding closed-loop system(s) in Step 2 is feasible.

One might wonder how complex to apply Theorem~\ref{th:comparison-vector-general} and Corollary~\ref{co:comparison-vector-general} in stability analysis since the conditions are quite complicated. In fact, the functions of $F(\cdot)$ and $G(\cdot)$ are normally chosen in such a way that they are quite similar in terms of structure. In most cases, they can be chosen to share the same structure. Then Condition 1 in Theorem~\ref{th:comparison-vector-general} and Corollary~\ref{co:comparison-vector-general} is, by default, true. Moreover, the function of $f(\cdot)$ is normally quite similar to $g(\cdot)$ in order to guarantee that Condition 2 in Theorem~\ref{th:comparison-vector-general} and Corollary~\ref{co:comparison-vector-general} can be satisfied. Conditions 3 and 5 in Theorem~\ref{th:comparison-vector-general} and Corollary~\ref{co:comparison-vector-general} are relatively easy to check given $F(\cdot)$ and $f(\cdot)$. Although Condition 4 in Theorem~\ref{th:comparison-vector-general} and Corollary~\ref{co:comparison-vector-general} is generally hard to check because the number of solutions to~\eqref{eq:equality-derivative} and~\eqref{eq:equality-Dini} could be infinite, the solutions normally share a common structure (see the following Section~\ref{sec:directed}). Then Condition 4 is also easy to analyze.

In the following part of this paper, we let $F(\cdot)$ and $G(\cdot)$ share the same structure. Thus, Condition 1 in Corollary~\ref{co:comparison-vector-general} is satisfied apparently. In addition, Conditions 3 and 5 in Corollary~\ref{co:comparison-vector-general} are satisfied under the properly designed $F(\cdot)$ and $f(\cdot)$. Therefore, we will not explicitly mention Conditions 1, 3, and 5 in the following of the paper when Corollary~\ref{co:comparison-vector-general} is used.

\section{Finite-time Consensus Under a Directed Switching Interaction Graph}\label{sec:directed}
In this section, we focus on the study of finite-time consensus for multi-agent networks with dynamics given in~\eqref{eq:kinmatics} under a directed switching interaction graph. In the following of this paper, we only focus on the case of $m=1$ (\emph{i.e.}, the one-dimensional case) for the simplicity of presentation. However, because condition~\eqref{eq:lipshitz} holds for each dimension, similar results can be obtained for $m>1$ (\emph{i.e.}, the high-dimensional case) by applying a similar analysis to each dimension.

Before moving on, we need the following lemmas.

\begin{lemma}\label{lem:negative-semi}
Let $n$ be any positive integer and $x_i\in\re,~i=1,\cdots,n$. There exists a positive number $q_n$ such that
\begin{equation}\label{eq:negative-semi}
-\sum_{i=1}^n x_i^2+\sum_{i=2}^{n}x_ix_{i-1}
\leq -q_n \sum_{i=1}^n x_i^2.
\end{equation}
\end{lemma}
\proof When $n=1$, the lemma holds apparently because
\begin{align*}
&-\sum_{i=1}^n x_i^2+\sum_{i=2}^{n}x_ix_{i-1}
=-x_1^2\leq -\sum_{i=1}^n x_i^2.
\end{align*}
We next consider the case when $n\geq 2$.

Since $-\sum_{i=1}^n x_i^2+\sum_{i=2}^{n}x_ix_{i-1}=0$ if $x_i=0$ for all $i=1,\cdots,n$, in order to prove the lemma, it is equivalent to show that $-\sum_{i=1}^n x_i^2+\sum_{i=2}^{n}x_ix_{j-1}$ is negative whenever there exists $i\in\{1,\cdots,n\}$ such that $x_i\neq 0$. Note that
\begin{align*}
   &-\sum_{i=1}^n x_i^2+\sum_{i=2}^{n}x_ix_{i-1}  =   -\frac{1}{2}\left[2\sum_{i=1}^n x_i^2-\sum_{i=2}^{n}2x_ix_{i-1}\right]\\
   =  & -\frac{1}{2} \Big[x_1^2+x_n^2+(x_2-x_1)^2+(x_n-x_{n-1})^2\\
   &~~~~~~~~ +2\sum_{i=2}^{n-1}(x_{i+1}-x_i)^2\Big].
\end{align*}
Therefore, $-\sum_{i=1}^n x_i^2+\sum_{i=2}^{n}x_ix_{j-1}\leq 0$. Note also that $-\sum_{i=1}^n x_i^2+\sum_{i=2}^{n}x_ix_{j-1}=0$ if and only if
\begin{align*}
  x_1=&0,\quad x_n=0,\\
  x_i-x_{i-1}=&0,~i=2,\cdots,n,
\end{align*}
which implies that $-\sum_{i=1}^n x_i^2+\sum_{i=2}^{n}x_ix_{j-1}=0$ if and only if $x_i=0,~i=1,\cdots,n$. Combining with the fact $-\sum_{i=1}^n x_i^2+\sum_{i=2}^{n}x_ix_{j-1}\leq 0$ shows that $-\sum_{i=1}^n x_i^2+\sum_{i=2}^{n}x_ix_{j-1}$ is negative whenever there exists $i\in\{1,\cdots,n\}$ such that $x_i\neq 0$. This completes the proof of the lemma.
\endproof

\begin{lemma}~\cite{XiaoWangJia08}\label{lem:XiaoWangJia}
Let $x_1,\cdots,x_p\geq 0$ and $\alpha^\star\in(0,1)$. Then $\sum_{i=1}^px_i^{\alpha^\star}\geq (\sum_{i=1}^p x_i)^{\alpha^\star}$.
\end{lemma}

\begin{lemma}~\cite{Khalil02}\label{lem:comparison-lemme}
Consider the scalar differential equation
\[\dot{\mu}=g(t,\mu),\quad \mu(t_0)=\mu_0,\]
where $g(t,\mu)$ is continuous in $t$ and locally Lipschitz in $\mu$ for all $t\geq 0$ and all $\mu\in J\in\re$. Let $[t_0,T)$ (T could be infinity) be the maximal interval of existence of the solution $\mu(t)$, and suppose $\mu(t)\in J$ for all $t\in[t_0, T)$. Let $\nu(t)$ be a continuous function whose upper Dini derivative $D^+\nu(t)$ satisfies the differential inequality
\[D^+\nu(t)\leq g(t,\nu(t)),\quad \nu(t_0)\leq \mu_0\]
with $\nu(t)\in J$ for all $t\in[t_0,T)$. Then $\nu(t)\leq \mu(t)$ for all $t\in[t_0,T)$.
\end{lemma}

With the aid of Corollary~\ref{co:comparison-vector-general} and Lemmas~\ref{lem:negative-semi},~\ref{lem:XiaoWangJia}, and~\ref{lem:comparison-lemme}, we next show that finite-time consensus is achieved for~\eqref{eq:kinmatics} using~\eqref{eq:control-new} in the following theorem.

\begin{theorem}\label{th:directed-finite-2}
Assume that the directed interaction graph $\Gcal_i,i=0,1,\ldots,$ has a directed spanning tree. Using~\eqref{eq:control-new} for~\eqref{eq:kinmatics}, $\abs{r_i(t)-r_j(t)}\to 0$ in finite time for any initial state if $\beta\geq\frac{\gamma+\epsilon_1}{\underline{a}q_{n-1}}+\epsilon_2$, where $q_n$ is the maximal positive number such that~\eqref{eq:negative-semi} holds and $\epsilon_k,~k=1,2$, is any positive constant.
\end{theorem}
\proof
The closed-loop system of~\eqref{eq:kinmatics} using~\eqref{eq:control-new} can be written as
\begin{align}\label{eq:clp-theorem}
\dot{r}_i=f(t,r_i)-\beta\sum_{j=1}^n a_{ij}(t) \text{sig}(r_i-r_j)^{\alpha(\abs{r_i-r_j})}.
\end{align}
Note that $\abs{r_i-r_j}\to 0$ in finite time if and only if $\max_{i}r_i-\min_{i}r_i\to 0$ in finite time. We next employ Corollary~\ref{co:comparison-vector-general} to prove the theorem. We first show that the finite-time convergence of~\eqref{eq:kinmatics} using~\eqref{eq:control-new} can be guaranteed by the finite-time convergence of another well-designed closed-loop system. Then we prove that the designed closed-loop system indeed converges in finite time.

\textit{The Use of Corollary~\ref{co:comparison-vector-general}:}

(1): Choose a proper nonnegative function $G(r)$ and compute $D^+G(r)$, where $r=[r_1,\cdots,r_n]^T$. Here we choose $G(r)\defeq \max_{i}r_i-\min_{i}r_i$. We first show that $G(r)=0$ for any $t\geq \underline{t}$ if $G(r)=0$ at $t=\underline{t}$. Whenever $G(r)=0$ at some time $\underline{t}$, it follows that $r_i(\underline{t})=r_j(\underline{t})$ for all $i,j=1,\cdots,n$. It thus follows from~\eqref{eq:clp-theorem} that $\dot{r}_i=\dot{r}_j,~\forall i,j=1,\cdots,n,$ at time $\underline{t}$. Consequently, $r_i(t)=r_j(t),~\forall i,j=1,\cdots,n$, for any $t\geq \underline{t}$. That is, $G(r)=0$ for any $t\geq \underline{t}$.

Based on the definition of upper Dini derivative, we have
\begin{align*}
D^+G(r)=&\limsup_{h\to 0^+} \frac{1}{h}\left\{G[r(t+h)]-G[r(t)]\right\}\\
=&\limsup_{h\to 0^+}
\frac{1}{h}\big\{[\max_i r_{i}(t+h)-\max_i r_{i}(t)]\\
&~~~~~~~~~~~~~~-[\min_ir_{i}(t+h)-\min_ir_{i}(t)]\big\}\\
=&D^+\max_i r_i-D^+\min_ir_i.
\end{align*}
Note that
\begin{align*}
&\limsup_{h\to 0^+}
\frac{1}{h}[\max_i r_{i}(t+h)-\max_i r_{i}(t)]\\
=&\lim_{\epsilon\to 0,\varepsilon\to 0^+}\sup\left\{\frac{\max_i r_{i}(t+h)-\max_i r_{i}(t)}{h}:h\in(\varepsilon-\epsilon,\varepsilon+\epsilon)\right\}\\
=&\lim_{\epsilon\to 0,\varepsilon\to 0^+}\left(\max_{i\in\arg\max_{\ell}{r}_\ell}\frac{r_{i}(t+h)- r_{i}(t)}{h}:h\in(\varepsilon-\epsilon,\varepsilon+\epsilon)\right)\\
=&\max_{i\in\arg\max_{\ell}{r}_\ell}\left(\lim_{\epsilon\to 0,~\varepsilon\to 0^+}\frac{r_{i}(t+h)- r_{i}(t)}{h}:h\in(\varepsilon-\epsilon,\varepsilon+\epsilon)\right).
\end{align*}
Because $h\to 0^+$ when $h\in(\varepsilon-\epsilon,\varepsilon+\epsilon)$, $\epsilon\to 0$, and $\varepsilon\to 0^+$, it follows that
\[\lim_{\epsilon\to 0,~\varepsilon\to 0^+,~h\in(\varepsilon-\epsilon,\varepsilon+\epsilon)}\frac{r_{i}(t+h)- r_{i}(t)}{h}=D^+ r_i(t).\]
Therefore,
$D^+\max_i r_i=\max_{i\in\arg\max_{\ell}{r}_\ell}D^+ r_i(t).$
Similarly,
$D^+\min_i r_i=\max_{i\in\arg\min_{\ell}{r}_\ell}D^+ r_i(t).$
Note that $\dot{r}_i(t)$ is continuous with respect to $t$ everywhere except for the case when the network topology is switching. Since the case only happens at some distinct time instants, $D^+ r_i(t)$ and $\dot{r}_i(t)$ are equivalent almost everywhere. That is, $D^+G(r)$ and $\max_{i\in\arg\max_{\ell}{r}_\ell}\dot{r}_i(t)-\max_{i\in\arg\min_{\ell}{r}_\ell}\dot{r}_i(t)$ are equivalent almost everywhere. From~\eqref{eq:clp-theorem}, we can obtain that
\begin{align*}
&\max_{i\in\arg\max_{\ell}{r}_\ell}\dot{r}_i(t)\\
= &\max_{i\in\arg\max_{\ell}{r}_\ell}\left[f(t,r_i)-\beta\sum_{j=1}^n a_{ij}(t) \text{sig}(r_i-r_j)^{\alpha(\abs{r_i-r_j})}\right]\\
\leq &\max_{i\in\arg\max_{\ell}{r}_\ell}f(t,r_i),
\end{align*}
where we have used the fact $\max_i r_i\geq r_j,~j=1,\cdots,n$, to derive the last inequality.
Similarly, it can be obtained that $\max_{i\in\arg\min_{\ell}{r}_\ell}\dot{r}_i(t)\geq \max_{i\in\arg\min_{\ell}{r}_\ell}f(t,r_i)$. When the directed interaction graph $\Gcal_i,i=0,1,\ldots,$ has a directed spanning tree, it can be obtained that $\{j|j\in\Ncal_i,i\in\arg\max_{\ell}r_\ell\}\bigcup\{j|j\in\Ncal_i,i\in\arg\min_{\ell}r_\ell\}\neq\varnothing$, where $\varnothing$ denotes the empty set.
We then rewrite $D^+G(r)$ based on the following two different cases:
\begin{itemize}
\item[(i)] When $\{j|j\in\Ncal_i,i\in\arg\max_{\ell}r_\ell\}\neq\varnothing$, by recalling that $D^+ r_i(t)$ and $\dot{r}_i(t)$ are equivalent almost everywhere and $\max_{i\in\arg\min_{\ell}{r}_\ell}\dot{r}_i(t)\geq \max_{i\in\arg\max_{\ell}{r}_\ell}f(t,r_i)$, it follows that $D^+G(r)\leq \max_{i\in\arg\max_{\ell}r_\ell}\dot{r}_i(t)-\max_{i\in\arg\min_{\ell}r_\ell}f(t,r_i)$ almost everywhere. Therefore, we have
$
D^+G(r)\leq \max_{i\in\arg\max_{\ell}r_\ell}\Bigg[f(t,r_i)-\beta\sum_{i=1}^n a_{ij}(t)\text{sig}(r_i-r_j)^{\alpha(\abs{r_i-r_j})}\Bigg]-\max_{i\in\arg\min_{\ell}r_\ell}f(t,r_i)< (\gamma+\epsilon_1)(\max_{i\in\arg\max_{\ell}r_\ell} r_i-\max_{i\in\arg\min_{\ell}r_\ell} r_i)-\max_{i\in\arg\max_{\ell}r_\ell} \beta\underline{a} \text{sig}(r_i-\max_{r_j< r_i}r_j)^{\alpha(\abs{r_i-\max_{r_j< r_i}r_j})}
$
holds almost everywhere when $G(r)\neq 0$, where~\eqref{eq:lipshitz} was used to derive the last inequality and $\epsilon_1$ is any positive constant.\newline
\item[(ii)] When $\{j|j\in\Ncal_i,i\in\arg\min_{\ell}r_\ell\}\neq \varnothing$ and $\{j|j\in\Ncal_i,i\in\arg\max_{\ell}r_\ell\}=\varnothing$, by following a similar analysis as that in Case (i), it can be obtained that $D^+G(r)\leq \max_{i\in\arg\max_{\ell}r_\ell}f(t,r_i)-\max_{i\in\arg\min_{\ell}r_\ell}\dot{r}_i(t)$ almost everywhere. Therefore, we have
$
D^+G(r)\leq \max_{i\in\arg\max_{\ell}r_\ell}f(t,r_i)-\max_{i\in\arg\min_{\ell}r_\ell}\Bigg[f(t,r_i)-\beta\sum_{i=1}^n a_{ij}(t)\text{sig}(r_i-r_j)^{\alpha(\abs{r_i-r_j})}\Bigg]
< (\gamma+\epsilon_1)(\max_{i\in\arg\max_{\ell}r_\ell}r_i-\max_{i\in\arg\min_{\ell}r_\ell}r_i)+\max_{i\in\arg\min_{\ell}r_\ell} \beta\underline{a}\text{sig}(r_i-\min_{r_j> r_i}r_j)^{\alpha(\abs{r_i-\min_{r_j> r_i}r_j})}
$
holds almost everywhere when $G(r)\neq 0$, where~\eqref{eq:lipshitz} was used again to derive the last inequality.
\end{itemize}
Note that $G(r)$ is not necessarily a nonincreasing function.

(2) Propose a proper function $F(\xi)$ satisfying Condition 1 and the corresponding closed-loop system(s) satisfying~\eqref{eq:equality-Dini} such that Conditions 2-5 are satisfied, where $\xi=[\xi_1,\ldots,\xi_n]^T$. Define $\hat{\gamma}\defeq \gamma+\epsilon_1$. Let $F(\xi)\defeq\max_i\xi_i-\min_i\xi_i$ and consider the closed-loop system given by
\begin{itemize}
\item[(i)] When $\{j|j\in\Ncal_i,i\in\arg\max_{\ell}\xi_\ell\}\neq\varnothing$,
\begin{align}\label{eq:Fxi-case1}
  &D^+F(\xi)
=\hat{\gamma}(\max_{i\in\arg\max_{\ell}\xi_\ell} \xi_i-\max_{i\in\arg\min_{\ell}\xi_\ell} \xi_i)\notag\\
&-\max_{i\in\arg\max_{\ell}\xi_\ell} \beta\underline{a} \text{sig}(\xi_i-\max_{\xi_j< \xi_i}\xi_j)^{\alpha(\abs{\xi_i-\max_{\xi_j< \xi_i}\xi_j})}
\end{align}
\item[(ii)] When $\{j|j\in\Ncal_i,i\in\arg\min_{\ell}\xi_\ell\}\neq \varnothing~\text{and}~\{j|j\in\Ncal_i,i\in\arg\max_{\ell}\xi_\ell\}=\varnothing$,
\begin{align}\label{eq:Fxi-case2}
&D^+F(\xi)= \hat{\gamma}(\max_{i\in\arg\max_{\ell}\xi_\ell}\xi_i-\max_{i\in\arg\min_{\ell}\xi_\ell}\xi_i)\notag\\
&+\max_{i\in\arg\min_{\ell}\xi_\ell} \beta\underline{a}\text{sig}(\xi_i-\min_{\xi_j> \xi_i}\xi_j)^{\alpha(\abs{\xi_i-\min_{\xi_j> \xi_i}\xi_j})}.
\end{align}
\end{itemize}
It then follows that Condition 2 in Corollary~\ref{co:comparison-vector-general} is satisfied. It remains unclear if Condition 4 in Corollary~\ref{co:comparison-vector-general} is satisfied.

By following a similar analysis to that of $D^+G(r)$ in Step 1, it can be obtained that \[D^+F(\xi)=\max_{i\in\arg\max_{\ell}{\xi}_\ell}D^+ \xi_i(t)-\max_{i\in\arg\min_{\ell}{\xi}_\ell}D^+ \xi_i(t).\]
Then~\eqref{eq:Fxi-case1} becomes
\begin{align}\label{eq:Fxi-case1-rewritten}
&\max_{i\in\arg\max_{\ell}{\xi}_\ell}\left[D^+ \xi_i(t)-\hat{\gamma}\xi_i+\beta\underline{a} \text{sig}(\xi_i-\max_{\xi_j< \xi_i}\xi_j)^{\alpha(\abs{\xi_i-\max_{\xi_j< \xi_i}\xi_j})}\right]\notag\\
&=\max_{i\in\arg\min_{\ell}{\xi}_\ell}\left[D^+ \xi_i(t)-\hat{\gamma}\xi_i\right].
\end{align}
Define
\[
\Gamma_i(t,\xi)\defeq D^+ \xi_i(t)-\hat{\gamma}\xi_i+\beta\underline{a} \text{sig}(\xi_i-\max_{\xi_j< \xi_i}\xi_j)^{\alpha(\abs{\xi_i-\max_{\xi_j< \xi_i}\xi_j})}.
\]
Then $\max_{i\in\arg\min_{\ell}{\xi}_\ell}\Gamma_i(t,\xi)=D^+ \xi_i(t)-\hat{\gamma}\xi_i$ because $\max_{\xi_j< \xi_i}\xi_j$ does not exist for $i\in\arg\min_{\ell}{\xi}_\ell$. Then~\eqref{eq:Fxi-case1-rewritten} becomes
\begin{align}\label{eq:Gamma_i-function}
\max_{i\in\arg\max_{\ell}{\xi}_\ell}\Gamma_i(t,\xi)=\max_{i\in\arg\min_{\ell}{\xi}_\ell}\Gamma_i(t,\xi).
\end{align}
In order to guarantee that~\eqref{eq:Gamma_i-function} holds for an arbitrary $\xi(0)$ and an arbitrary interaction graph, $\Gamma_i(t,\xi)$ should have the same structure of $\xi$ for all $i$ because otherwise $\max_{i\in\arg\max_{\ell}{\xi}_\ell}\Gamma_i(t,\xi)-\max_{i\in\arg\min_{\ell}{\xi}_\ell}\Gamma_i(t,\xi)$ cannot always be constant, which then results in a contradiction. As a consequence, $\Gamma_i(t,\xi)$ should satisfy
\begin{align}\label{eq:rho}
\Gamma_i(t,\xi)=\rho(t,\xi),\quad i=1,\cdots,n,
\end{align}
for some function $\rho(t,\xi)$ such that the (Caratheodory) solution to~\eqref{eq:Fxi-case1} exists. Recalling the definition of $\Gamma_i(t,\xi)$, it then follows from~\eqref{eq:rho} that
\begin{align}\label{eq:solutions}
D^+ \xi_i(t)=&\hat{\gamma}\xi_i-\beta\underline{a} \text{sig}(\xi_i-\max_{\xi_j< \xi_i}\xi_j)^{\alpha(\abs{\xi_i-\max_{\xi_j< \xi_i}\xi_j})}\notag\\
&+\rho(t,\xi),\quad i=1,\cdots,n.
\end{align}
In order to guarantee that the (Caratheodory) solution to~\eqref{eq:Fxi-case1} exists, it is required that $\rho(t,\xi)$ be integrable. That is, $\int_{0}^t\rho(\tau,\xi(\tau))\text{d}\tau$ is well-defined. Let's further define $\hat{\xi}_i(t)\defeq\xi_i(t)+\lambda$, where $\lambda$ is the (Caratheodory) solution to the following equation
\begin{equation}\label{eq:D+lambda}
D^+ \lambda=\hat{\gamma}\lambda-\rho(t,\xi),
\end{equation}
Here the (Caratheodory) solution to~\eqref{eq:D+lambda} exists because $\int_{0}^t\rho(\tau,\xi(\tau))\text{d}\tau$ is well-defined.
Then we have
\begin{align}\label{eq:xi_i_new}
D^+ \hat{\xi}_i(t)=&D^+\xi_i(t)+D^+\lambda\notag\\
=&\hat{\gamma}\xi_i-\beta\underline{a} \text{sig}(\xi_i-\max_{\xi_j< \xi_i}\xi_j)^{\alpha(\abs{\xi_i-\max_{\xi_j< \xi_i}\xi_j})}\notag\\
&+\rho(t,\xi)+\hat{\gamma}\lambda-\rho(t,\xi)\notag\\
=&\hat{\gamma}\xi_i-\beta\underline{a} \text{sig}(\xi_i-\max_{\xi_j< \xi_i}\xi_j)^{\alpha(\abs{\xi_i-\max_{\xi_j< \xi_i}\xi_j})}\notag\\
&+\hat{\gamma}[\hat{\xi}_i(t)-\xi_i(t)]\notag\\
=&\hat{\gamma}\hat{\xi}_i-\beta\underline{a} \text{sig}(\hat{\xi}_i-\max_{\hat{\xi}_j< \hat{\xi}_i}\hat{\xi}_j)^{\alpha(\abs{\hat{\xi}_i-\max_{\hat{\xi}_j< \hat{\xi}_i}\hat{\xi}_j})},\notag\\
&~~~~~~~~~~~ i=1,\cdots,n,
\end{align}
where we have used the fact that $\hat{\xi}_i(t)-\xi_i(t),~i=1,\cdots,n$, are identical to derive the last equality. Since the set of the discontinuity points for the right-hand side of~\eqref{eq:xi_i_new} has measure zero, the solution to~\eqref{eq:xi_i_new} in the Caratheodory sense is the same as the solution to
\begin{align}\label{eq:clp-gen-cmp-directed}
\dot{\hat{\xi}}_i(t)=&\hat{\gamma}\hat{\xi}_i-\beta\underline{a} \text{sig}(\hat{\xi}_i-\max_{\hat{\xi}_j< \hat{\xi}_i}\hat{\xi}_j)^{\alpha(\abs{\hat{\xi}_i-\max_{\hat{\xi}_j< \hat{\xi}_i}\hat{\xi}_j})},\notag\\
&~~~~~~~~~~~~~~~ i=1,\cdots,n,
\end{align}
in the Caratheodory sense. Clearly, in the Caratheodory sense, the solution to~\eqref{eq:Fxi-case1} exists since the solution to~\eqref{eq:clp-gen-cmp-directed} exists. Define $F(\hat{\xi})\defeq\max_i\hat{\xi}_i-\min_i\hat{\xi}_i$. $F(\xi)$ under~\eqref{eq:Fxi-case1} and $F(\hat{\xi})$ under~\eqref{eq:clp-gen-cmp-directed} are always equal.

Similarly, when $\{j|j\in\Ncal_i,i\in\arg\min_{\ell}\xi_\ell\}\neq \varnothing~\text{and}~\{j|j\in\Ncal_i,i\in\arg\max_{\ell}\xi_\ell\}=\varnothing$, the (Caratheodory) solution to~\eqref{eq:Fxi-case2} exists. In addition, $F(\xi)$ under~\eqref{eq:Fxi-case2} and $F(\hat{\xi})$ under
\begin{align}\label{eq:clp-gen-cmp-directed-2}
  \dot{\hat{\xi}}_i=\hat{\gamma} \hat{\xi}_i+\beta\underline{a}\text{sig}(\hat{\xi}_i-\min_{\hat{\xi}_j> \hat{\xi}_i}\hat{\xi}_j)^{\alpha(\abs{\hat{\xi}_i-\min_{\hat{\xi}_j> \hat{\xi}_i}\hat{\xi}_j})}
\end{align}
are always equal.

By combining the previous arguments, $F(\xi)$ goes to zero in finite time for any solution $\xi$ to the closed-loop system switching between~\eqref{eq:Fxi-case1} and~\eqref{eq:Fxi-case2} if and only if $F(\hat{\xi})$ goes to zero in finite time for the solution $\hat{\xi}$ to the closed-loop system with the same switching when~\eqref{eq:Fxi-case1} and~\eqref{eq:Fxi-case2} are replaced by, respectively,~\eqref{eq:clp-gen-cmp-directed} and~\eqref{eq:clp-gen-cmp-directed-2}. Therefore, in order to guarantee that Condition 4 in Corollary~\ref{co:comparison-vector-general} is satisfied, a sufficient condition is that $\max_i \hat{\xi}_i-\min_i\hat{\xi}_i\to 0$ in finite time for the system with an \emph{arbitrary}
switching between~\eqref{eq:clp-gen-cmp-directed} and~\eqref{eq:clp-gen-cmp-directed-2}, which will be shown next.

\textit{The proof of $\max_i\hat{\xi}_i(t)-\min_i\hat{\xi}_i(t)\to 0$ in finite time for the system with an arbitrary switching between~\eqref{eq:clp-gen-cmp-directed} and~\eqref{eq:clp-gen-cmp-directed-2}:}



The proof can be divided into two steps:
\begin{itemize}
\item[1] Prove that $\max_i\hat{\xi}_i(t)-\min_i\hat{\xi}_i(t)\to 0$ as $t\to\infty$;
\item[2] Prove that $\max_i\hat{\xi}_i(t)-\min_i\hat{\xi}_i(t)\to 0$ in finite time.
\end{itemize}
Here the proof of Step 2 depends on the statement in Step 1.

\textit{Step 1: Prove that $\max_i\hat{\xi}_i(t)-\min_i\hat{\xi}_i(t)\to 0$ as $t\to\infty$.} Before moving on, we first analyze an important property for the system with an arbitrary switching between~\eqref{eq:clp-gen-cmp-directed} and~\eqref{eq:clp-gen-cmp-directed-2}. That is, if $\hat{\xi}_i(t_0)\leq\hat{\xi}_j(t_0)$, then $\hat{\xi}_i(t)\leq\hat{\xi}_j(t)$ for any $t>t_0$. This property is due to the fact that whenever $\hat{\xi}_i(t_0)=\hat{\xi}_j(t_0)$ at some $t_0$, $\dot{\hat{\xi}}_i(t_0)=\dot{\hat{\xi}}_j(t_0)$ and thus $\hat{\xi}_i(t)=\hat{\xi}_j(t)$ for any $t>t_0$. This property plays an important role in the following proof. Without loss of generality, we label the agents in such a way that $\hat{\xi}_1(0)\leq \hat{\xi}_2(0)\leq\ldots\leq \hat{\xi}_n(0)$ throughout the following proof. Combining with the property implies that $\hat{\xi}_1(t)\leq \hat{\xi}_2(t)\leq\ldots\leq \hat{\xi}_n(t)$ for any $t>0$.

Consider the Lyapunov function candidate given by
\begin{align}\label{eq:Lyapunov-G}
\tilde{G}(\hat{\xi})\defeq\sum_{i=1}^{n-1}\int_0^{\hat{\xi}_{i+1}-\hat{\xi}_i} \tau^{\alpha(\abs{\tau})}\text{d}\tau.
\end{align}
Noting that $\tau^{\alpha(\abs{\tau})}$ is continuous with respect to $\tau$, it follows that $\int_0^{\hat{\xi}_{i+1}-\hat{\xi}_i} \tau^{\alpha(\abs{\tau})}\text{d}\tau$ is differentiable with respect to $\hat{\xi}_{i+1}-\hat{\xi}_i$. Therefore, $\tilde{G}(\hat{\xi})$ is regular. Then $\tilde{G}(\hat{\xi})$ is a nonpathological function~\cite{BacciottiCeragioli06}. In the sense of Caratheodory solutions for the system with an arbitrary switching between~\eqref{eq:clp-gen-cmp-directed} and~\eqref{eq:clp-gen-cmp-directed-2}, if the nonpathological derivative of $\tilde{G}(\hat{\xi})$ is negative definite, then the system with an arbitrary switching between~\eqref{eq:clp-gen-cmp-directed} and~\eqref{eq:clp-gen-cmp-directed-2} is globally asymptotically stable (c.f. Lemma 1 in~\cite{BacciottiCeragioli06}). Because the nonpathological derivative of a nonpathological function is essentially equivalent to the set-valued derivative of the same function (see Definition $4$ in~\cite{BacciottiCeragioli06}), we next use the notation of the set-valued derivative in the proof.

The set-valued derivative of $\tilde{G}(\hat{\xi})$ is given by
$
L_F\tilde{G}(\hat{\xi})=
K\left[\sum_{j=1}^{n-1}(\hat{\xi}_{i+1}-\hat{\xi}_i)^{\alpha(\abs{\hat{\xi}_{i+1}-\hat{\xi}_{i}})}(\dot{\hat{\xi}}_{i+1}-\dot{\hat{\xi}}_{i})\right],
$
where $K[\cdot]$ is the differential inclusion~\cite{Filippov88}. Before analyzing $L_F\tilde{G}(\hat{\xi})$, we first analyze $\dot{\hat{\xi}}_{i+1}-\dot{\hat{\xi}}_{i}$ under~\eqref{eq:clp-gen-cmp-directed} and~\eqref{eq:clp-gen-cmp-directed-2}. When~\eqref{eq:clp-gen-cmp-directed} is satisfied, we have
$
\dot{\hat{\xi}}_{i+1}-\dot{\hat{\xi}}_{i}
= \hat{\gamma}(\hat{\xi}_{i+1}-\hat{\xi}_i)-\beta\underline{a}\text{sig}(\hat{\xi}_{i+1}-\max_{\hat{\xi}_j< \hat{\xi}_{i+1}}\hat{\xi}_j)^{\alpha(\abs{\hat{\xi}_{i+1}-\max_{\hat{\xi}_j< \hat{\xi}_{i+1}}\hat{\xi}_j})}
+\beta\underline{a}\text{sig}(\hat{\xi}_i-\max_{\hat{\xi}_j< \hat{\xi}_i}\hat{\xi}_j)^{\alpha(\abs{\hat{\xi}_i-\max_{\hat{\xi}_j< \hat{\xi}_i}\hat{\xi}_j})}.
$
Because $\hat{\xi}_{i+1}\geq \hat{\xi}_i$ (see the first paragraph of the step), we have
$
\abs{\hat{\xi}_{i+1}-\max_{\hat{\xi}_j< \hat{\xi}_{i+1}}\hat{\xi}_j}
=\hat{\xi}_{i+1}-\max_{\hat{\xi}_j< \hat{\xi}_{i+1}}\hat{\xi}_j
\geq \hat{\xi}_{i+1}-\hat{\xi}_i\geq 0.
$
It thus follows that
$
\hat{\xi}_{i+1}-\hat{\xi}_i\leq (\hat{\xi}_{i+1}-\hat{\xi}_i)^{\alpha(\abs{\hat{\xi}_{i+1}-\hat{\xi}_i})}.
$
When $\hat{\xi}_1(t)< \hat{\xi}_2(t)<\ldots< \hat{\xi}_n(t)$, it follows that $\max_{\hat{\xi}_j< \hat{\xi}_{i+1}}\hat{\xi}_j=\hat{\xi}_{i}$. Define $\bar{\delta}_i\defeq (\hat{\xi}_{i+1}-\hat{\xi}_i)^{\alpha(\abs{\hat{\xi}_{i+1}-\hat{\xi}_i})}$. Then we obtain
\begin{align*}
  &\sum_{j=1}^{n-1}(\hat{\xi}_{i+1}-\hat{\xi}_i)^{\alpha(\abs{\hat{\xi}_{i+1}-\hat{\xi}_i})}(\dot{\hat{\xi}}_{i+1}-\dot{\hat{\xi}}_{i})\\
  =&\sum_{j=1}^{n-1}(\hat{\xi}_{i+1}-\hat{\xi}_i)^{\alpha(\abs{\hat{\xi}_{i+1}-\hat{\xi}_i})}\Bigg[\hat{\gamma}(\hat{\xi}_{i+1}-\hat{\xi}_i)\\
  &-\beta\underline{a}\text{sig}(\hat{\xi}_{i+1}-\hat{\xi}_i)^{\alpha(\abs{\hat{\xi}_{i+1}-\hat{\xi}_i})}
  +\beta\underline{a}\text{sig}(\hat{\xi}_i-\hat{\xi}_{i-1})^{\alpha(\abs{\hat{\xi}_i-\hat{\xi}_{i-1}})}\Bigg]\\
  \leq &\sum_{j=1}^{n-1}(\hat{\xi}_{i+1}-\hat{\xi}_i)^{\alpha(\abs{\hat{\xi}_{i+1}-\hat{\xi}_i})}\Bigg[\hat{\gamma}(\hat{\xi}_{i+1}-\hat{\xi}_i)^{\alpha(\abs{\hat{\xi}_{i+1}-\hat{\xi}_i})}\\
  &-\beta\underline{a}\text{sig}(\hat{\xi}_{i+1}-\hat{\xi}_i)^{\alpha(\abs{\hat{\xi}_{i+1}-\hat{\xi}_i})}
  +\beta\underline{a}\text{sig}(\hat{\xi}_i-\hat{\xi}_{i-1})^{\alpha(\abs{\hat{\xi}_i-\hat{\xi}_{i-1}})}\Bigg]\\
=& -(\beta\underline{a}-\hat{\gamma})\sum_{i=1}^{n-1} \bar{\delta}_i^{2}+\beta\underline{a}\sum_{i=2}^{n-1}\bar{\delta}_i\bar{\delta}_{i-1}\\
  \leq & -(\beta\underline{a}q_{n-1}-\hat{\gamma})\sum_{i=1}^{n-1} \bar{\delta}_i^{2}
  \leq -\underline{a}q_{n-1}\epsilon_2 \sum_{i=1}^{n-1} \bar{\delta}_i^{2},
\end{align*}
where we have used Lemma~\ref{lem:negative-semi} and $\beta\geq\frac{\hat{\gamma}}{\underline{a}q_{n-1}}+\epsilon_2$ to derive the last two inequalities. When $\hat{\xi}_1(t)< \hat{\xi}_2(t)<\ldots< \hat{\xi}_n(t)$ does not hold and $\tilde{G}(\hat{\xi})\neq 0$ as time $t$, we can always find a set of indices such that $\hat{\xi}_{k_1}(t)< \hat{\xi}_{k_2}(t)<\ldots< \hat{\xi}_{k_\ell}(t)$ and $\hat{\xi}_i(t)\in\{\hat{\xi}_{k_1}(t),\cdots,\hat{\xi}_{k_\ell}(t)\}$ for any $i\notin\{k_1,\cdots,k_\ell\}$. Then the Lyapunov function candidate~\eqref{eq:Lyapunov-G} can be equivalently written as
\[\tilde{G}(\hat{\xi})\defeq\sum_{i=1}^{\ell-1}\int_0^{\hat{\xi}_{k_{i+1}}-\hat{\xi}_{k_i}} \tau^{\alpha(\abs{\tau})}\text{d}\tau.\]
Define $\tilde{\delta}_i\defeq (\hat{\xi}_{k_{i+1}}-\hat{\xi}_{k_{i}})^{\alpha(\abs{\hat{\xi}_{k_{i+1}}-\hat{\xi}_{k_i}})}$. Similarly, it can be obtained that
$
\sum_{j=1}^{\ell-1}(\hat{\xi}_{k_{i+1}}-\hat{\xi}_{k_i})^{\alpha(\abs{\hat{\xi}_{k_{i+1}}-\hat{\xi}_{k_i}})}(\dot{\hat{\xi}}_{k_{i+1}}-\dot{\hat{\xi}}_{k_{i}})
\leq -\underline{a}q_{n-1}\epsilon_2 \sum_{i=1}^{\ell-1} \tilde{\delta}_i^{2}.
$
Since $\bar{\delta}_i=0$ if and only if $\hat{\xi}_{i+1}-\hat{\xi}_i=0$, $\sum_{i=1}^{n-1}\bar{\delta}_i^2$ remains unchanged if those agents with the same state are considered one agent. This implies that $\sum_{i=1}^{\ell-1} \tilde{\delta}_i^{2}=\sum_{i=1}^{n-1} \bar{\delta}_i^{2}$.
Therefore, $\max L_F \tilde{G}(\hat{\xi})$ is always negative definite when~\eqref{eq:clp-gen-cmp-directed} is satisfied.
When~\eqref{eq:clp-gen-cmp-directed-2} is satisfied and $\hat{\xi}_1(t)< \hat{\xi}_2(t)<\ldots< \hat{\xi}_n(t)$ holds (respectively, $\hat{\xi}_1(t)< \hat{\xi}_2(t)<\ldots< \hat{\xi}_n(t)$ does not hold),
by following a similar analysis to that of the case when~\eqref{eq:clp-gen-cmp-directed} is satisfied and $\hat{\xi}_1(t)< \hat{\xi}_2(t)<\ldots< \hat{\xi}_n(t)$ (respectively, $\hat{\xi}_1(t)< \hat{\xi}_2(t)<\ldots< \hat{\xi}_n(t)$ does not hold), we can obtain that $\max L_F \tilde{G}(\hat{\xi})\leq -\underline{a}q_{n-1}\epsilon_2 \sum_{i=1}^{n-1} \bar{\delta}_i^{2}$, which implies that $\max L_F \tilde{G}(\hat{\xi})$ is also negative definite.


From the previous analysis, we know that $\max L_F \tilde{G}(\hat{\xi})$ is always negative definite for the system with an arbitrary switching between~\eqref{eq:clp-gen-cmp-directed} and~\eqref{eq:clp-gen-cmp-directed-2}. It then follows from Lemma $1$ in~\cite{BacciottiCeragioli06} that $\hat{\xi}_i(t)-\hat{\xi}_j(t)\to 0,~\forall i,j\in\{1,\cdots,n\},$ as $t\to\infty$. Equivalently, $\max_i\hat{\xi}_i(t)-\min_i\hat{\xi}_i(t)\to 0$ as $t\to\infty$.

\textit{Step 2: Prove that $\max_i\hat{\xi}_i(t)-\min_i\hat{\xi}_i(t)\to 0$ in finite time.} Because $\max_i\hat{\xi}_i(t)-\min_i\hat{\xi}_i(t)\to 0$ as $t\to\infty$ (shown in Step 1), it follows that there exists a time instant $\overline{t}$ such that $\max_i\hat{\xi}_i(t)-\min_i\hat{\xi}_i(t)<1$ for any $t\geq \overline{t}$, which implies that $\abs{\xi_i(t)-\xi_j(t)}<1$ for any $t\geq \overline{t}$. Based on~\eqref{eq:alpha-func}, for $t\geq \overline{t}$,~\eqref{eq:clp-gen-cmp-directed} and~\eqref{eq:clp-gen-cmp-directed-2} can be rewritten as
\begin{align}\label{eq:clp-gen-cmp-directed-new}
  \dot{\hat{\xi}}_i=\hat{\gamma} \hat{\xi}_i-\beta\underline{a} \text{sig}(\hat{\xi}_i-\max_{\hat{\xi}_j< \hat{\xi}_i}\hat{\xi}_j)^{\alpha^\star}
\end{align}
and
\begin{align}\label{eq:clp-gen-cmp-directed-2-new}
  \dot{\hat{\xi}}_i=\hat{\gamma} \hat{\xi}_i+\beta\underline{a}\text{sig}(\hat{\xi}_i-\min_{\hat{\xi}_j> \hat{\xi}_i}\hat{\xi}_j)^{\alpha^\star}.
\end{align}
The Lyapunov function candidate~\eqref{eq:Lyapunov-G} can be rewritten as
\begin{equation}\label{eq:G-new}
\tilde{G}(\hat{\xi})=\sum_{i=1}^n\int_{0}^{\hat{\xi}_{i+1}-\hat{\xi}_i} \tau^{\alpha^\star}\text{d}\tau.
\end{equation}
By following a similar analysis to that of $L_F\tilde{G}(\hat{\xi})$, we can obtain that
\begin{align*}
\max L_F\tilde{G}(\hat{\xi})\leq &-\underline{a}q_{n-1}\epsilon_2\sum_{i=1}^{n-1} (\hat{\xi}_{i+1}-\hat{\xi}_i)^{2\alpha^\star}\\
=&-\underline{a}q_{n-1}\epsilon_2\sum_{i=1}^{n-1} [(\hat{\xi}_{i+1}-\hat{\xi}_i)^{1+\alpha^\star}]^{\frac{2\alpha^\star}{1+\alpha^\star}}\\
\leq &-\underline{a}q_{n-1}\epsilon_2[\sum_{i=1}^{n-1} (\hat{\xi}_{i+1}-\hat{\xi}_i)^{1+\alpha^\star}]^{\frac{2\alpha^\star}{1+\alpha^\star}},
\end{align*}
where we have used Lemma~\ref{lem:XiaoWangJia} to derive the last inequality.
Noting that $\tilde{G}(\hat{\xi})=\frac{1}{(1+\alpha^\star)}\sum_{i=1}^n (\hat{\xi}_{i+1}-\hat{\xi}_i)^{1+\alpha^\star},$ it follows that
\begin{equation}\label{eq:maxLF}
\max L_F\tilde{G}(\hat{\xi})\leq-\underline{a}q_{n-1}\epsilon_2(1+\alpha^\star)^{\frac{2\alpha^\star}{1+\alpha^\star}}[\tilde{G}(\hat{\xi})]^{\frac{2\alpha^\star}{1+\alpha^\star}}.
\end{equation}
Let's write $\tilde{G}(\hat{\xi})$ as $\tilde{G}(t)$ for simplicity. For $t\geq \overline{t}$, although $\dot{\tilde{G}}(t)$ is discontinuous at some time instants, it is always integrable because $\tilde{G}(t)\in[0,\tilde{G}(\overline{t})]$ is bounded and the set of the discontinuity points has measure zero~\cite{Apostol74}. It then follows from~\eqref{eq:maxLF} that
\begin{align*}
\tilde{G}(t+h)-\tilde{G}(t)=&\int_{t}^{t+h} \dot{\tilde{G}}(\tau)\text{d}\tau\\
\leq & -h\underline{a}q_{n-1}\epsilon_2(1+\alpha^\star)^{\frac{2\alpha^\star}{1+\alpha^\star}}\min_{\tau\in[t,t+h]}\tilde{G}(t).
\end{align*}
It then follows from the definition of upper Dini derivative (c.f. Section~\ref{sec: pre}) that
\begin{align*}
D^+ \tilde{G}(t) \leq -\underline{a}q_{n-1}\epsilon_2(1+\alpha^\star)^{\frac{2\alpha^\star}{1+\alpha^\star}}[\tilde{G}(t)]^{\frac{2\alpha^\star}{1+\alpha^\star}}.
\end{align*}
When $t\geq \overline{t}$, it follows from Lemma~\ref{lem:comparison-lemme} that $\tilde{G}(t)$ is upper bounded by $\mu(t)$ satisfying
\[\dot{\mu}(t)=-\underline{a}q_{n-1}\epsilon_2(1+\alpha^\star)^{\frac{2\alpha^\star}{1+\alpha^\star}}[\mu(t)]^{\frac{2\alpha^\star}{1+\alpha^\star}},\quad \mu(\overline{t})=\tilde{G}(\overline{t}).\]
Because $\frac{2\alpha^\star}{1+\alpha^\star}\in(0,1)$, by computation, we have
\begin{align*}
&\frac{1+\alpha^\star}{1-\alpha^\star}[\mu(t)]^{\frac{1-\alpha^\star}{1+\alpha^\star}}\\
=&\frac{1+\alpha^\star}{1-\alpha^\star}[\mu(\overline{t})]^{\frac{1-\alpha^\star}{1+\alpha^\star}}-\underline{a}q_{n-1}\epsilon_2(1+\alpha^\star)^{\frac{2\alpha^\star}{1+\alpha^\star}}(t-\overline{t}).
\end{align*}
Therefore, $\mu(t)\to 0$ in finite time. Because $\tilde{G}(t)$ is nonnegative and $\tilde{G}(t)$ is upper bounded by $\mu(t)$, $\tilde{G}(t)\to 0$ in finite time. Equivalently, $\max_i\hat{\xi}_i(t)-\min_i\hat{\xi}_i(t)\to 0$ in finite time.

Combining all previous arguments completes the proof.
\endproof

\begin{myremark}
Although the solutions to~\eqref{eq:Fxi-case1} in Step 2 of the proof of Theorem~\ref{th:directed-finite-2} are not unique, they share some nice common features [\textit{i.e.,} they all satisfy~\eqref{eq:solutions}]. Similarly, the solutions to~\eqref{eq:Fxi-case2} in Step 2 of the proof of Theorem~\ref{th:directed-finite-2} also share some similar nice common features. Ultimately, the study of finite-time consensus for the system switching between~\eqref{eq:Fxi-case1} and~\eqref{eq:Fxi-case2} is converted to the study of finite-time consensus for the system switching between~\eqref{eq:clp-gen-cmp-directed} and~\eqref{eq:clp-gen-cmp-directed-2}. The example shows that the verification of Condition 4 in Theorem~\ref{th:comparison-vector-general} and Corollary~\ref{co:comparison-vector-general} is not as complex as it appears to be.
\end{myremark}

\begin{myremark}
In the proof of Theorem~\ref{th:directed-finite-2}, the stability and finite-time convergence of~\eqref{eq:kinmatics} using~\eqref{eq:control-new} under a (general) directed switching interaction graph is shown to be guaranteed by the stability and finite-time convergence of the closed-loop system switching between~\eqref{eq:clp-gen-cmp-directed} and~\eqref{eq:clp-gen-cmp-directed-2}. The interaction graph associated with~\eqref{eq:clp-gen-cmp-directed} [respectively,~\eqref{eq:clp-gen-cmp-directed-2}] is constructed in such a way that each agent has at most one neighbor whose state is the maximum of those states that are smaller than its own state (respectively, the minimum of those states that are larger than its own state).
\end{myremark}

It can be observed that the closed-loop system of~\eqref{eq:kinmatics} using~\eqref{eq:control-new} is nonlinear and the stability analysis is, in general, difficult because the closed-loop system is nonlinear and switching in the presence of unknown terms. By using Corollary~\ref{co:comparison-vector-general}, the stability of~\eqref{eq:kinmatics} using~\eqref{eq:control-new} can be guaranteed by the stability of another nonlinear system, whose stability can be analyzed. In addition, the unknown dynamics do not appear in the new nonlinear system.


When the unknown inherent nonlinear dynamics do not exist, we have the following corollary.
\begin{corollary}\label{co:conv-finite-directed}
Consider agents with single-integrator kinematics
\begin{equation}
\dot{x}_i=u_i,\quad i=1,\cdots,n\label{eq:1st-order-system}
\end{equation}
where $x_i\in\re$ is the state of the $\ith$ agent and $u_i\in\re$ is the control input for the $\ith$ agent.
Let the consensus algorithm for~\eqref{eq:1st-order-system} be given by
\begin{equation}\label{finite-switching-11}
u_i =-\epsilon \sum_{j=1}^n a_{ij}(t) \text{sig}(x_i-x_j)^{\alpha^\star},
\end{equation}
where $\epsilon$ is any positive constant, $0<\alpha^\star<1$, and $a_{ij}(t)$ is the $(i,j)$th entry of the adjacency matrix $\Acal(t)$ at time $t$. Assume that the interaction graph $\Gcal_i,i=0,1,\ldots,$ has a directed spanning tree. Then $x_i(t)-x_j(t)\to 0$ in finite time.
\end{corollary}
\proof The corollary is a direct result of Theorem~\ref{th:directed-finite-2}.
\endproof

\begin{myremark}\label{re:Jointnotwork}
In~\cite{WangXiao10}, it is shown that finite-time consensus is reached if the directed \emph{fixed} graph has a directed spanning tree and each strongly connected component is \emph{detail-balanced}. In Corollary~\ref{co:conv-finite-directed}, we show that finite-time consensus is reached if the directed \emph{switching} graph has a directed spanning tree at each time interval. Then it is natural to ask if the condition on the interaction graph can be further relaxed to the case when the interaction graph has a directed spanning tree in some joint fashion in terms of a union of its time-varying graph topologies (see, e.g.,~\cite{RenBeard05_TAC} for details). Unfortunately, the relaxed condition on the interaction graph containing a directed spanning tree in some joint fashion, in general, cannot guarantee finite-time consensus, or even asymptotical consensus. For example, consider three agents with $r_1(0)=0,~r_2(0)=1,$ and $r_3=2$. Let $a_{23}=1$ for a period of $t_1$ such that $r_2(t_1)=r_3(t_1)$ with other entries in $\Acal$ equal to zero. Then let $a_{31}=1$ for a period of $t_2$ such that $r_3(t_1+t_2)=r_1(t_1+t_2)$ with other entries in $\Acal$ equal to zero. By continuing a similar process, both $\max_i r_i$ and $\min_i r_i$ keep unchanged, which implies that consensus cannot be achieved even if the interaction graph has a directed spanning tree jointly.
\end{myremark}

\begin{myremark}
Another interesting finite-time consensus algorithm for~\eqref{eq:1st-order-system} is given in~\cite{XiaoWCG09} by
\begin{equation}\label{eq:XiaoWCG}
u_i=-\epsilon \text{sig}\left[\sum_{j=1}^n a_{ij}(t) (x_i-x_j)\right]^\alpha,
\end{equation}
where $\epsilon$ is a positive constant. It is shown in~\cite{XiaoWCG09} that finite-time consensus is reached for~\eqref{eq:1st-order-system} using~\eqref{eq:XiaoWCG} when the directed \emph{fixed} interaction graph $\Gcal$ has a directed spanning tree. When the interaction graph $\Gcal_i,i=0,1,\ldots,$ has a directed spanning tree, by following a similar analysis to that in the proof of Theorem~\ref{th:directed-finite-2}, it can be shown that consensus is reached in finite time for~\eqref{eq:1st-order-system} using~\eqref{eq:XiaoWCG}. Note again that the relaxed condition on the interaction graph containing a directed spanning tree in some joint fashion as mentioned in Remark~\ref{re:Jointnotwork}, in general, cannot guarantee finite-time consensus.
\end{myremark}



Inspired by the main idea behind the proof of Theorem~\ref{th:directed-finite-2} for~\eqref{eq:control-new}, another finite-time consensus algorithm for~\eqref{eq:kinmatics} is given by
\begin{equation}\label{eq:control-cont}
u_i=-k\sum_{j=1}^n a_{ij}(t)(r_i-r_j)-\beta\sum_{j=1}^n a_{ij}(t)\text{sig}(r_i-r_j)^{\alpha^\star},
\end{equation}
where $k$ and $\beta$ are positive constants, and $\alpha^\star\in(0,1)$ is a positive constant. By following a similar analysis to that in the proof of Theorem~\ref{th:directed-finite-2}, it can be shown that consensus is reached in finite time when the interaction graph $\Gcal_i,i=0,1,\ldots,$ has a directed spanning tree, $k>\frac{\gamma}{\underline{a}q_{n-1}}$, and $\beta>0$. The main idea behind~\eqref{eq:control-cont} is that the linear term $-k\sum_{j=1}^n a_{ij}(t)(r_i-r_j)$ is used to compensate for the unknown inherent nonlinear dynamics such that asymptotic consensus can be achieved while the nonlinear term $-\beta\sum_{j=1}^n a_{ij}(t)\text{sig}(r_i-r_j)^{\alpha^\star}$ is used to guarantee finite-time consensus in the absence of the unknown inherent nonlinear dynamics. In the absence of the unknown inherent nonlinear dynamics, the use of the linear term is helpful (or at least harmless) for the nonlinear term to achieve finite-time consensus. Similarly, in the presence of the unknown inherent nonlinear dynamics, the use of the nonlinear term is helpful (or at least harmless) for the linear term to achieve asymptotic consensus. More precisely, both terms are harmless to each other in terms of each individual objective. Here an interesting idea is provided regarding the design of proper algorithms motivated by those algorithms under some similar situations. It is possible that a summation of various algorithms can be used to achieve some combinatorial objective if those algorithms are chosen properly. In addition, the stability analysis tool proposed in Section~\ref{sec:General_comparison_lemma} provides important insights as to the connection among the closed-loop systems using those algorithms.


We next consider a special case when the nonlinear part of~\eqref{eq:control-cont}, namely $-\beta\sum_{j=1}^n a_{ij}(t)\text{sig}(r_i-r_j)^{\alpha^\star}$, does not exist.
\begin{corollary}\label{co:asymp-consensus-nonlinear}
Assume that the directed interaction graph $\Gcal_i,i=0,1,\ldots,$ has a directed spanning tree. Using
\[
u_i=-k\sum_{j=1}^n a_{ij}(t)(r_i-r_j)
\]
for~\eqref{eq:kinmatics}, $\abs{r_i(t)-r_j(t)}\to 0$ as $t\to\infty$ if $k\geq\frac{\gamma+\epsilon_1}{\underline{a}q_{n-1}}+\epsilon_2$, where $q_n$ is the maximal positive number such that~\eqref{eq:negative-semi} holds and $\epsilon_k,~k=1,2$, is any positive constant.
\end{corollary}
\proof The proof is similar to that of Theorem~\ref{th:directed-finite-2}.
\endproof
\begin{myremark}
In contrast to~\cite{YuChenCao11} where asymptotical consensus for~\eqref{eq:kinmatics} was studied under some fixed interaction graphs, Corollary~\ref{co:asymp-consensus-nonlinear} presents more general results for the case of a directed switching interaction graph.
\end{myremark}

\section{Simulation}\label{sec:sim}
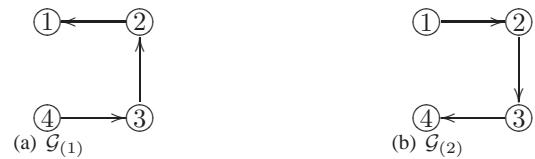
\begin{figure}[hhhhtb]
   \begin{center}
   \begin{tabular}{cc}
      \subfigure[$\Gcal_{(1)}$]{
         \begin{xy}*!D\xybox{\xymatrix{
            &*+[o][F-]{1}
            & *+[o][F-]{2}\ar[l]
              \\
            & *+[o][F-]{4}\ar[r]
            & *+[o][F-]{3}\ar[u]
         }}\end{xy}}
         \hspace{2cm}
      \subfigure[$\Gcal_{(2)}$]{
        \begin{xy}*!D\xybox{\xymatrix{
            &*+[o][F-]{1}\ar[r]
            & *+[o][F-]{2}\ar[d]
              \\
            & *+[o][F-]{4}
            & *+[o][F-]{3}\ar[l]
        }}\end{xy}}
     \end{tabular}
    \caption{Directed graphs $\Gcal_{(1)}$ and $\Gcal_{(2)}$. Both of them have a directed spanning tree. An arrow from $j$ to $i$ denotes that agent $j$ is a neighbor of agent $i$.}
    \label{fig:L-st-switch-directed}
   \end{center}
\end{figure}

In this section, a simulation example is presented to validate the theoretical result presented
in Section~\ref{sec:directed}. We consider a group of $4$ agents in the one-dimensional space (\textit{i.e.,} $m=1$) with the interaction graph switching from $\{\Gcal_{(1)},\Gcal_{(2)}\}$ (see Fig.~\ref{fig:L-st-switch-directed}) every $0.5$ seconds. By computation, $q_3=0.6910$. In particular, $a_{ij}(t)=1$ if $(j,i)\in\mathcal{W}(t)$ and $a_{ij}(t)=0$ otherwise. By choosing $f(t,r_i)=\sin(r_i)$, it then follows that $\abs{f(t,r_i)-f(t,r_j)}\leq \abs{r_i-r_j}$. We then choose $\beta=3$, which satisfies the condition on $\beta$ given in Theorem~\ref{th:directed-finite-2}. We further choose $\alpha^\star=0.8$ and $r(0)=[\frac{\pi}{2},-\frac{\pi}{2},-\frac{\pi}{2},-\frac{\pi}{2}]^T$.

\begin{figure}
\begin{center}
\includegraphics[width=.5\textwidth]{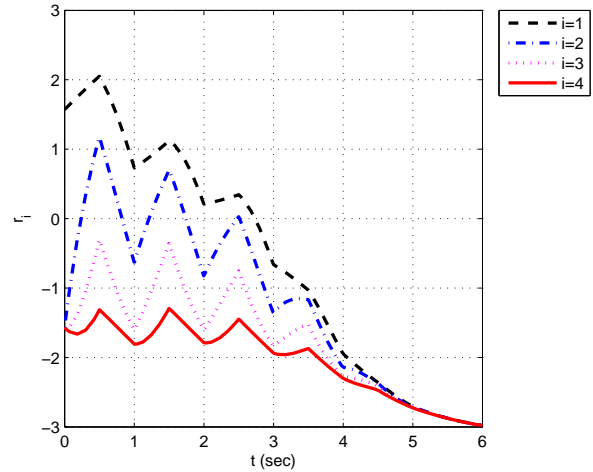}    
\caption{The trajectories of the four agents using~\eqref{eq:control-new}.}  
\label{fig:smooth-dir}                                 
\end{center}                                 
\end{figure}

\begin{figure}
\begin{center}
\includegraphics[width=.5\textwidth]{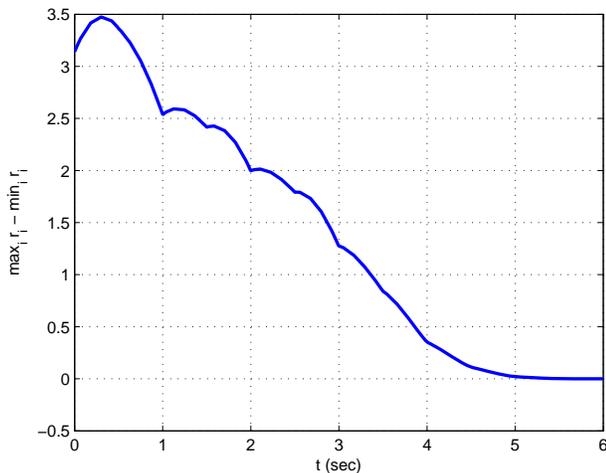}    
\caption{The evolution of $\max_i r_i-\min_ir_i$ with respect to $t$.}  
\label{fig:smooth-dir-diff}                                 
\end{center}                                 
\end{figure}

\begin{figure}
\begin{center}
\includegraphics[width=.5\textwidth]{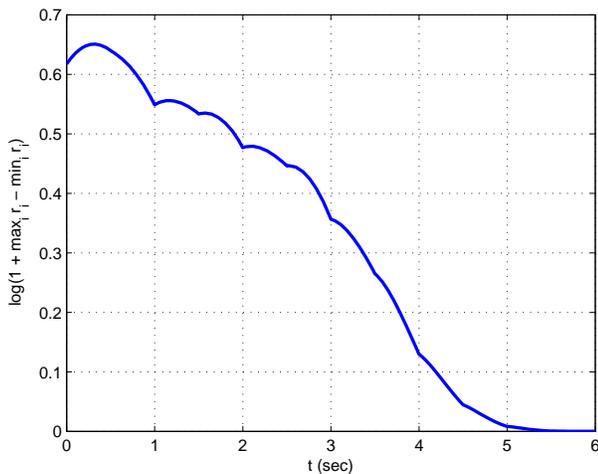}    
\caption{The evolution of $\log (1+\max_i r_i-\min_ir_i)$ with respect to $t$.}  
\label{fig:log-diff}                                 
\end{center}                                 
\end{figure}

Fig.~\ref{fig:smooth-dir} shows the trajectories of the four agents under the proposed consensus algorithm~\eqref{eq:control-new}. The evolution of $G(r)=\max_i r_i-\min_ir_i$ is given in Fig.~\ref{fig:smooth-dir-diff}. To better show the finite-time convergence, the evolution of $\log (1+\max_i r_i-\min_ir_i)$ is shown in Fig.~\ref{fig:log-diff}. Recall that $\max_i r_i-\min_ir_i=0$ if and only if $\log (1+\max_i r_i-\min_ir_i)=0$. It can be noticed from Figs.~\ref{fig:smooth-dir-diff} and~\ref{fig:log-diff} that all agents reach consensus in
finite time. Moreover, it can be seen from Fig.~\ref{fig:smooth-dir-diff} that $G(r)$ might increase at some time intervals, indicating that $G(r)$ is not a Lyapunov function.

%
%

\section{Conclusion}\label{sec:conclusion}
In this paper, we studied finite-time consensus of multi-agent networks with unknown inherent nonlinear dynamics. First, we proposed a novel stability tool based on a generalized comparison
lemma. With the aid of the novel stability analysis tool, we analyzed the stability of the closed-loop system using the proposed distributed nonlinear consensus algorithm by comparing the original closed-loop system with some well-designed closed-loop system that can guarantee finite-time consensus. In particular, the stability and finite-time convergence using the proposed nonlinear consensus algorithm under a (general) directed switching interaction graph were shown to be guaranteed by those of some special well-designed nonlinear closed-loop system under some special directed switching interaction graph. As a byproduct, in the absence of the unknown inherent nonlinear dynamics, the proposed nonlinear consensus algorithm can still guarantee finite-time consensus under a directed switching interaction graph.

\bibliographystyle{IEEEtran}
\bibliography{refs}

\end{document}